\documentclass[12pt]{article}
\usepackage{graphicx, amsfonts, amsmath,  amsthm, amssymb, latexsym, url}
\usepackage{array}
\usepackage{psboxit}

\newcommand\comment[1]{}                
\renewcommand\comment[1]{\textrm{[#1]}}           

\newtheorem{theorem}{Theorem}

\newtheorem{lemma}{Lemma}[section]

\newtheorem{proposition}[theorem]{Proposition}

\theoremstyle{remark}

\theoremstyle{definition}
\newtheorem{definition}[theorem]{Definition}

\renewcommand\emptyset{\varnothing}
\renewcommand{\epsilon}{\varepsilon}
\renewcommand{\phi}{\varphi}

\renewcommand\tilde{\widetilde}

\def\det{\mathop{\rm det}}

\newcommand\Aff{\operatorname{Aff}} \renewcommand\dim{\operatorname{\mathsf{dim}}} 
\renewcommand\ker{\operatorname{\mathfrak{Ker}}}

 \newcommand\conv{\operatorname{\mathsf{conv}}} 
 \newcommand\aff{\operatorname{\mathsf{aff}}} 
\newcommand\lin{\operatorname{\mathsf{lin}}}  
\newcommand\vertex{\operatorname{\mathsf{vert}}}  
  
 \newcommand\rank{\operatorname{\mathsf{rank}}} 
  
 \newcommand\Quad{\operatorname{\mathtt{QForm}}} 
\newcommand\Quot{\operatorname{Quot}} \newcommand\Sym{\operatorname{\mathsf{Sym}}} 
\newcommand\QuadInv{\operatorname{\mathsf{QuadInv}}}

\newcommand\Spec{\operatorname{\mathsf{Spec}}}

\newcommand\cER{\mathcal{ER}} 
\newcommand\QP{\mathsf{QP}}

\newcommand\0{\mathbf{0}}
\newcommand\aaa{\mathbf{a}}

\newcommand\cc{\mathbf{c}}
\newcommand\e{\mathbf{e}}
\newcommand\bd{\mathbf{d}}

\newcommand\m{\mathbf{m}}

\newcommand\oo{\mathbf{o}}  
\newcommand\p{\mathbf{p}}

\newcommand\rr{\mathbf{r}}

\newcommand\ttt{\mathbf{t}}
\newcommand\vv{\mathbf{v}}  

\newcommand\x{\mathbf{x}}
\newcommand\y{\mathbf{y}}
\newcommand\z{\mathbf{z}}

\newcommand\cE{\mathcal{E}}
\newcommand\cF{\mathcal{F}}

\newcommand\cQ{\mathcal{Q}}

\newcommand\cS{\mathcal{S}}

\newcommand\cV{\mathsf{V}}

\newcommand\bE{\textsf{\slshape E}}

\newcommand\R{\mathbb{R}}     
\newcommand\Q{\mathbb{Q}}     
\newcommand\Z{\mathbb{Z}}     
\newcommand\E{\mathbb{E}}     

\begin{document}

\author{Mathieu Dutour, Robert  Erdahl, and Konstantin Rybnikov}
\title{Perfect Delaunay Polytopes in Low Dimensions}
\date{\today}
\maketitle

\begin{abstract}
A lattice Delaunay polytope is known as perfect if the only ellipsoid, that can be 
circumscribed about it, is its Delaunay sphere. Perfect Delaunay polytopes are in 
one-to-one correspondence with arithmetic equivalence classes of positive quadratic 
functions on $\Z^n$, that can be recovered in a unique way from its minimum over $\Z^n$ 
and all of its representations. We develop a structural theory of such polytopes and 
describe all known perfect Delaunay polytopes in dimensions one through eight. We 
suspect that this list is complete. 

\end{abstract}

\section{\protect \large Introduction}\label{sec:introduction}
A point lattice is a discrete set of points in $\R^n$ such that the difference vectors form a subgroup of $\R^n$.
If $\Lambda$ is a point lattice in $\R^n$ ($n \ge 0$), then a convex polytope $P \subset \R^n$ 
is called a lattice polytope (or $\Lambda$-polytope)  if all of its vertices are lattice 
points. Consider the lattice $\mathbb{Z}^{n} \subset \R^n$, and a convex $\mathbb{Z}^{n}$-polytope $P$.  If 
$P$ can be circumscribed by an ellipsoid $\mathcal{E}=\{\x \in \R^n \; \vline\; \cQ_{\cE}(\x-\cc_{\cE})\le 
\rho_{\cE}^2 \}$  with no interior $\mathbb{Z}^{n}$-elements so that the  boundary 
$\mathbb{Z}^{n}$-elements of $\cE$ are exactly the vertices of $P$, we will say that $P$ is a 
\textit{Delaunay} polytope with respect to the form $\cQ_{\cE}$ defined by 
$\mathcal{E}$; more informally, we will say that a lattice polytope is \textit{Delaunay} if it 
can be circumscribed by an \textit{empty ellipsoid}.  \emph{Ellipsoid} is commonly used to 
refer to hypersurfaces, defined by positive definite quadratic forms, as well as solid bodies 
bounded by such surfaces:  the meaning of our usage will be clear from the context.     
Typically, there is a family of empty ellipsoids that can be circumscribed about a given 
Delaunay polytope $P$, but, if there is only one, so that $\mathcal{E}$ is uniquely determined 
by $P$, we will say that $P$ is a \textit{perfect Delaunay polytope} in $\Z^n$. Perfect 
Delaunay polytopes are also sometimes referred to as \emph{extreme}. Perfect Delaunay polytopes 
are fascinating geometrical objects -- examples are the six- and seven-dimensional Gosset 
(1900) polytopes  with $27$ and $56$ vertices respectively, which appear in the Delaunay 
tilings of the root lattices $E_{6}$ and $E_{7}$ (see e.g. Coxeter (1988) for a description). 
In this paper we describe, up to an isometry and a dilation, all known perfect Delaunay 
polytopes in $\R^n$ for $n \le 8$; we also present a study of the geometry and combinatorics of 
these polytopes. We suspect that the list of perfect Delaunay polytopes that we give here is 
complete for $n \le 8$. Erdahl (1975, 1992) proved that $0$ and $[0,1]$ are the only perfect 
Delaunay polytopes for $n\le5$; Dutour (2002) proved that there is only one perfect Delaunay 
polytope in $\R^6$ -- the Gosset $6$-polytope (Coxeter's $2_{21} $), which is described in 
Section \ref{sec:list:6}. Only two perfect Delaunay polytopes are known in $\R^7$: they are 
Gosset's $7$-polytope (Coxeter's $3_{21}$) and a $35$-tope found by Erdahl and Rybnikov (2002), 
which are described in Section \ref{sec:list7}. We list 27 8-dimensional perfect Delaunay 
polytopes, they are identified by numbers, 1 through 27; Section \ref{sec:list8} contains a 
detailed description of these polytopes.  There are infinite series of perfect Delaunay polytopes -- the first such series was found by 
Erdahl and Rybnikov in 2001 (Rybnikov, 2001). This series was further generalized in  (Erdahl, Ordine, Rybnikov 2004) to a 3-parametric series (where one
parameter is the dimension) of perfect Delaunay polytopes. Another infinite series has been found by Dutour (2005). Prior to 2001 only sporadic 
examples of perfect Delaunay polytopes had been known, besides the cases for $n\le 7$ mentioned above,
 all of them were found by  Deza, Grishukhin and Laurent (1992) and all of their examples were constructed as sections of the Leech and Barnes-Wall lattices.

\section{\protect \large Definitions and Notation}\label{sec:defs}
Formally speaking, the subject of this paper is the study of 0-sets of positive 
$\Q$-valued quadratic functions on free $\Z$-modules of finite rank. Since any free 
$\Z$-module of finite rank $\Lambda$  can be realized as a discrete subgroup of $\R^n$ for 
any $n \ge \rank \Lambda$, we can think of $\Lambda$ geometrically as a discrete set 
of vectors in vector space $\R^n$ or as a discrete set of points in affine space $\R^n$. Whenever we approach a $\Z$-module from this point of 
view, we call it a \textit{lattice}. 

\begin{definition}
A quadratic lattice is a pair $(\Lambda,\R)$, where $\Lambda$ is a 
free $\Z$-module of finite rank and 
$\cQ:\Lambda \rightarrow \R$ is a quadratic 
form.
\end{definition}

A function $F$ on a module is called quadratic if it can be written as $\Quad F +A$, 
where $\Quad F$ is a quadratic form and $A$ is an affine function. In general, we denote the 
quadratic form part of a polynomial $P$  by $\Quad P$. 

\begin{definition} An \emph{affine} quadratic lattice is a pair $\Aff(\Lambda,Q)$, where 
$\Lambda$ is free $\Z$-module of finite rank and
$Q:\Lambda \rightarrow \bE$ is a quadratic function valued in an extension $\bE$ of $\Q$.
\end{definition}

Sometimes, to stress that we consider a quadratic lattice, as opposed to an affine quadratic 
lattice we call the former linear or homogeneous quadratic lattice. $\Aff(\Lambda,Q)$ (resp.  
$(\Lambda,\cQ)$) is called degenerate if $\rank \Quad Q<\rank \Lambda$  (resp.  $\rank \cQ 
<\rank \Lambda$).  Quadratic lattices $(\Lambda,\cQ)$ and $(\Lambda_1,\cQ_1)$ are called 
isomorphic if there is a $\Z$-module isomorphism $L:\Lambda \rightarrow \Lambda_1$ such that 
$\cQ(\x)=\cQ_1(L\x)$ for any $\x \in \Lambda$. Affine quadratic lattices $\Aff(\Lambda, Q)$ and 
$\Aff(\Lambda_1, Q_1)$ are called isomorphic if there is a $\Z$-module isomorphism $L:\Lambda 
\rightarrow \Lambda_1$ and some $\z \in \Lambda_1$ such that $Q(\x)=Q_1(L\x-\z)$ for any $\x 
\in \Lambda$.  We allow quadratic forms and functions to take values not only in $\Z$, the 
ground ring of the modules under consideration, but also in $\Q$ and its extensions, unlike the 
classical case of (linear) quadratic modules, where forms are supposed to be valued in the 
ground ring of the module (as in, e.g.,  Serre 1973).  Since we restrict ourselves to   $\rank 
\Lambda < \infty$, any  $\Q$-valued function on $\Lambda$ can be rescaled into a $\Z$-valued 
quadratic function. 
 
Suppose  $Q_1$ and $Q_2$ are valued in $\R$. Two  affine quadratic lattices 
$\Aff(\Lambda_1,Q_1)$ and $\Aff(\Lambda_2,Q_2)$ are called equivalent up to scaling if there  
is a $\Z$-module isomorphism $L:\Lambda_1 \rightarrow \Lambda_2$, some $\z \in \Lambda_2$ and a 
real $c >0$ such that $Q_1(\x)=c\,Q_2(L\x-\z)$.

A quadratic form $\cQ$ induces a symmetric bilinear form on $\Lambda \times \Lambda$ by 
\[(\x,\y) \mapsto \frac{1}{2}\left\{\cQ[\x+\y]-\cQ[\x]-\cQ[\y]\right\};\] \noindent we will 
denote this bilinear form also by $\cQ$ -- normally, there is no confusion, since a quadratic 
form has one arguments, while the  corresponding bilinear form has two. We call a number  $a$ 
positive (resp. negative) if $a \ge 0$ (resp. $a \le 0$); same terminology is applied to 
functions. Thus, a form $\cQ$ is called positive if $\cQ[\x] \ge 0$ for any $\x$. A form  $\cQ$ 
is called positive definite if $\cQ[\x] > 0$ for any $\x\neq \0$. 

Suppose a quadratic function $Q$ is valued in $\R$.  A number $b \in \R$ is called the arithmetic minimum of an 
affine quadratic lattice $\Aff(\Lambda,Q)$ if $b=\underset{\z \in \Lambda}{\min}\: Q(\z)$. The vectors of 
$\Lambda$ on which the minimum of  $\Aff(\Lambda,Q)$ is attained are called the minimal vectors of  $\Aff(\Lambda,Q)$.
 The definition of the arithmetic minimum of an \emph{affine } quadratic lattice is slightly 
 different from  that of a homogeneous quadratic lattice: in the case of a 
 homogeneous quadratic lattice the minimum is taken over all non-zero vectors $\z \in \Lambda$. 

\begin{definition} Let $Q:\Lambda \rightarrow \R$ be a fixed quadratic function with positive $\Quad Q$,
  and let $X$ be a 
quadratic function with unknown coefficients. Let $b=\underset{\z \in \Lambda}{\min}\: Q(\z)$. 
The affine quadratic lattice $\Aff(\Lambda,Q)$ is called \emph{perfect} if the 
system of equations 
\[\{X(\m)=b\;\;\vline\;\; \m\: \mbox{is a minimal vector of } Q \}\]
has the only solution  $X=Q$. \end{definition}

When we refer to a function as perfect without specifying $\Lambda$,  the meaning of 
$\Lambda$ is clear from the context. If a concrete formula for the function is given, 
$\Lambda$ is presumed to be  $\Z^n$. We will often use a shorthand such as e.g. 
\emph{$k$-lattice} (polytope, cell, etc.), instead of\emph{ $k$-dimensional lattice} 
(polytope, cell, etc.). 

From now on all quadratic functions are valued in $\Q$ or $\R$. In this context a symmetric 
bilinear form is a scalar product on $\Lambda \otimes \Q \cong \Q^n$ or $\Lambda \otimes \R 
\cong \R^n$.  There are two canonical ways to describe an affine quadratic lattice, one by 
fixing the lattice to be $\Z^n$ and the other by fixing the quadratic form part of the function 
to be, say $\sum x_i^2$. The first method is more flexible, since it allows for quadratic forms 
of arbitrary signature.  Furthermore, in any kind of machine computations it is far more 
convenient to deal with the former representation. 

In each dimension there are only finitely many non-isomorphic perfect affine quadratic lattices, up to  scale. This follows from Voronoi's
L-type reduction theory (see e.g. Deza and Laurent, 1997). Namely, in each dimension there is a strict L-type reduction domain $\mathcal{D}$,  which has finitely many extreme rays. The quadratic part of a perfect quadratic function is always arithmetically equivalent to a form lying on an extreme ray of $\mathcal{D}$ (but not vice versa). This implies the finiteness.

\begin{proposition}
For $n=0$ the only perfect affine quadratic lattice is $(\Z^0,0)$. For $n=1$ the only 
perfect affine quadratic lattice, up to isomorphisms and scaling, is 
$\Aff(\Z^1,(x-\frac{1}{2})^2)$. 
\end{proposition}

Perfect quadratic functions are inhomogeneous analogs of perfect quadratic forms introduced in 
the middle of 19-th century by Korkine and Zolotareff (1873) and later studied by Voronoi, 
Barnes, Conway and Sloane, Stacey, Martinet, and others (see Martinet (2003) for a survey). The 
interest to perfect forms has been mostly fueled by the theorem, proven by Korkine and 
Zolotareff (1873), that forms that are extreme points of the ball packing density function are 
perfect. We prefer not to use the term inhomogeneous form, employed in some number-theoretic 
literature (e.g. Gruber, Lekkerkerker 1987) since a\emph{ form} is by definition a 
\emph{homogeneous polynomial. } 

\subsection{\protect \bf Delaunay Tilings} The language of Delaunay tilings provides a 
geometric way of thinking about quadratic functions with positive quadratic part. We denote the 
vertex set of a polytope $P$ by $\vertex P$.   A convex $\Lambda$-polytope $D$ is called a 
Delaunay polytope in a (linear) quadratic lattice $(\Lambda,\cQ)$, where $\cQ[\x] > 0$ for $\x 
\neq 0$, if there is an ellipsoid $Q(\x) \le 0$, with $\Quad Q=\cQ$,  whose boundary contains  
$\vertex D$, but no points of $\Lambda\,\diagdown\,\vertex D$. If $\dim D < \rank \Lambda$ such 
an ellipsoid is not unique; however, the intersection of any such ellipsoid with the affine 
span of $D$ is unique (for fixed $(\Lambda,\cQ)$). In particular, such an ellipsoid is unique 
when $D$ is of maximal dimension -- in this case this ellipsoid is called the Delaunay 
ellipsoid (or empty ellipsoid) of $D$.

For any $S \subset \Lambda$,  we denote the affine span of $S$ in $\Lambda \otimes 
\R$ by $\aff S$. The affine span of   $S$ in $\Lambda \otimes \Q$ is denoted by  
$\aff_{\Q} S$, and the lattice spanned by all vectors $\x-\y$, where $\x,\y \in S$, 
by $\aff_{\Z} S$. Note that when $\0 \in S$, the linear span of $S$, $\lin S$, and the affine span of $S$, $\aff S$, are the same.  Often lattices arise as sections of other lattices by affine 
subspaces of the ambient affine Euclidean space. $\Gamma \subset \R^n$ is called an 
affine lattice if  $\overrightarrow{\Gamma}=\{\x-\y \;\vline\; \x,\y \in \Gamma\}$ is a 
$\Z$-module of finite rank. In such situations it is convenient to have the notion of 
 isomorphism for affine lattices. Let $\Gamma \subset \R^n$ and  $\Gamma' \subset \R^m$ be 
affine lattices. A map $f:\Gamma \rightarrow \Gamma'$ is called  an affine isomorphism if  
there are $\oo \in \Gamma$ and $\oo' \in \Gamma'$ such that   $\x-\oo \mapsto f\x-\oo'$ is a 
$\Z$-module isomorphism from $\overrightarrow{\Gamma}$ onto $\overrightarrow{\Gamma'}$.  Two functions 
$\phi:\Gamma \rightarrow S$, $\psi:\Gamma' \rightarrow S$ on affine lattices $\Gamma$ and  
$\Gamma'$ are called arithmetically equivalent if there is an affine isomorphism $f:\Gamma 
\rightarrow \Gamma'$ such that $\phi(\x)=\psi(f\x)$. A $\Lambda$-polyhedron $P$  can be thought 
of as the indicator function, which is 1 on $P\cap \Lambda$ and 0 elsewhere on $\Lambda$. Then, 
arithmetic equivalence of lattice polyhedra is a special case of the arithmetic equivalence 
between functions.  If $\Gamma=\Gamma'=\Z^n$, the arithmetic equivalence is the same as the 
equivalence with respect to $Aff(n,\Z)$, the group of all transformations of type $\z \mapsto 
L\z+\ttt$, where $L \in GL(n,\Z)$ and $\ttt \in \Z^n$.

 It is a theorem of Delaunay (1924) that for a positive definite quadratic form 
 $\cQ:\Lambda \longrightarrow \R$ the space $\aff \Lambda$ is partitioned 
 into the relative interiors of Delaunay $\Lambda$-polytopes with respect to $\cQ$; this 
partition is organized so that the intersection of any family of Delaunay polytopes is again a 
Delaunay polytope (we add the empty polytope $\emptyset$ to the partition) -- in other words 
the resulting \textit{Delaunay tiling} is \textit{face-to-face.}  This theorem also says that a 
Delaunay tiling for $(\Lambda,\cQ)$ is unique.  In studying Delaunay tilings and L-types of 
lattices it is often beneficial not to restrict to positive definite forms, but to use the 
concept of Delaunay tiling with respect to any positive $\Q$-valued quadratic form. Since 
traditionally, in the geometric context, quadratic forms are valued  in $\R$, we will say a few 
words about the case where $\cQ$ is $\R$-valued.

\begin{definition} The rational closure of the cone of positive-definite quadratic forms on 
$\Lambda$ is the set of all positive $\R$-valued forms $\cQ$ that satisfy the condition  $\rank 
(\Lambda \cap \ker_{\R} \cQ)=\dim \ker_{\R} \cQ$. \end{definition} It is easy to see that the 
rational closure of the cone of positive-definite quadratic forms on $\Lambda$ is a convex cone 
over $\R$. When $\Lambda=\Z^n$ we denote by $\Sym(n,\R)$ the space of $\R$-valued quadratic 
forms on $\Z^n$, and by ${\Sym}_+(n,\R)$  cone of all positive-definite quadratic forms in 
$\Sym(n,\R)$.  Then the real closure  of ${\Sym}_+(n,\R)$ in $\Sym(n,\R)$ is denoted by 
$\overline{\Sym}_+(n,\R)$, and the rational closure of ${\Sym}_+(n,\R)$  by 
$\overline{\Sym}_+^{\Q}(n,\R)$. When we consider an arbitrary lattice $\Lambda$, rather than $\Z^n$, 
we write $\Sym(\Lambda,*)$ instead of $\Sym(n,*)$.

$\overline{\Sym}_+^{\Q}(n,\R)$ can also be described as the real cone spanned by rank-one forms in 
indeterminates $(x_1,\ldots,x_n)=\x$ of type $(\vv \cdot \x)^2$ where $\vv$ runs over  $\Z^n$ 
(see e.g. Dutour, Schuermann, Vallentin, 2006).

Sometimes the condition $\rank ( \Z^n \cap \ker_{\R} \cQ )=\dim \ker_{\R} \cQ$ is phrased as 
that $\cQ$ has rational kernel, although this expression can be somewhat misleading, since 
$\ker_{\R} \cQ \cap \Z^n$ is always a rational subspace of $\Q^n$.  Since, 
$\overline{\Sym}^{\Q}_+(n,\R) \cap \Sym(n,\Q)$ consists of all positive  forms with rational 
coefficients, perhaps, it would be more elegant to consider only $\Q$-valued forms on $\Q^n$, 
but we decided to follow the tradition and embed the cone of positive $\Q$-valued  forms  into 
$\overline{\Sym}^{\Q}_+(n,\R)$. Let us denote by $\tilde{\QP}_{+}(n,\R)$ the cone of all real 
quadratic polynomials on $\R^n$ whose quadratic form parts belong to 
$\overline{\Sym}^{\Q}_+(n,\R)$.
 It is easy to see that $\tilde{\QP}_{+}(n,\R)$ is a convex cone  in the space ${\QP}(n,\R)$ of quadratic polynomials on $\R^n$.

When $\cQ:\Lambda \rightarrow \R$  is positive semidefinite, 
it defines a tiling of $\Lambda \otimes \R$ only when 
$\rank (\ker_{\R} \cQ \cap \Lambda)=\dim \ker_{\R} \cQ$, i.e. 
when $\cQ \in \overline{\Sym}^{\Q}_+(\Lambda,\R)$ 
(see e.g. Dutour, Schuermann, Vallentin, 2006). In this case $\vertex P$ 
should be interpreted as $\partial P \cap  \Lambda$ rather than as the set of vertices of 
$P$ in the sense of geometry of $\R^n$; \
for simplicity, we still call elements of $\vertex P$ vertices.

\begin{definition} Let $(\Lambda,\cQ)$ be a quadratic lattice with  $\cQ \in 
\overline{\Sym}^{\Q}_+(\Lambda,\R)$.  A convex polyhedron $P \subset \Lambda \otimes \R$ is a 
Delaunay polyhedron for   $(\Lambda,\cQ)$  if there is a quadratic polynomial $E_P$ on $\Lambda 
\otimes \R$, with $\Quad E_P =\cQ$, such that $E_P(\z)=0$ for any $\z \in \Lambda \cap P$ and 
$E_P(\z)>0$ for any  $\z \in \Lambda \diagdown P$. \end{definition} 

In particular, the empty polytope $\emptyset$ and the whole space $\Lambda \otimes \R$ are 
Delaunay polyhedra in  $(\Lambda,0)$: polynomials $E_{\emptyset}=1$ and    $E_{\Lambda \otimes 
\R}=0$ can serve as witnesses. It is not difficult to show that when $\cQ:\Lambda 
\longrightarrow \Q$  is positive semidefinite, $\aff \Lambda$ is covered by  Delaunay 
\emph{polyhedra} of various dimensions: some of these polyhedra are polytopes and some are 
cylinders over Delaunay  polytopes of lower dimensions. We often refer to a Delaunay polyhedron 
as  a Delaunay cell. In the semidefinite case the relative interiors of Delaunay cells also  
form a face-to-face partition of $\Lambda \otimes \R$, but  the elements of $\Lambda$ can no 
longer be considered as $0$-cells of the tiling -- the tiling in this case does not have any 
$0$-cells unless $\rank \Lambda =0$.    We denote the set of all cells of the Delaunay tiling 
of $\Lambda \otimes \R$ with respect to a semidefinite form $\cQ$ by $Del(\Lambda,\cQ)$.  
$Del(\Lambda,\cQ)$ has a poset structure, namely, $F \preceq C$ if and only if $F \subset 
\partial C$.  Furthermore, since $\emptyset,\:\Lambda \otimes \R \in Del(\Lambda,\cQ)$, it is a 
lattice. In discussions of concrete Delaunay polytopes   it is often more convenient to refer 
to faces by their vertex sets.  The partial order on  $Del(\Lambda,\cQ)$ induces a partial 
order on the vertex  (in the generalized sense explained above)  sets of Delaunay cells of 
$Del(\Lambda,\cQ)$;  we denote the resulting poset by $\cV Del(\Lambda,\cQ)$.  We will need the 
notion of Delaunay tiling for degenerate quadratic lattices only in Subsection 
\ref{subsec:layer}, so, for the exception of that part of the paper, the reader may safely 
assume that $\cQ$ is positive definite and all Delaunay cells are polytopes.  

For formal definitions and detailed information on Delaunay tilings of lattices we refer to 
Deza and Laurent (1997). We only remark that the  Delaunay tilings for lattices are classically 
(Delaunay, 1924) defined with the Euclidean norm $x_1^2+\ldots+x_n^2$ (in geometry of numbers 
the norm of a vector is its squared length), but are most effectively studied by isomorphicly 
mapping the lattice $\Lambda$ onto $\mathbb{Z}^{n}$, and replacing the Euclidean norm by a 
positive quadratic form $\cQ$ that makes $(\Z^n,\cQ)$ and $(\Lambda,\sum x_i^2)$ equivalent. 
This allows us to think in terms of Euclidean lattices, i.e. geometrically, but compute in 
terms of quadratic forms.

\begin{definition}
 Let $(\Lambda,\cQ)$ be a quadratic lattice with  $\cQ \in \overline{\Sym}^{\Q}_+(\Lambda,\R)$. Suppose $P  \in Del(\Lambda,\cQ)$.  $P$ is called perfect if its Delaunay ellipsoid (or elliptic cylinder, if $\rank \cQ < \rank \Lambda$) with respect to 
 $\cQ$ is the only quadric circumscribed about $P$ in $\Lambda \times \R$. 
\end{definition}  

Indeed, the notion of perfection, that was introduced in 19th century by the Italian school of 
algebraic geometry,  is\textit{ independent of the Delaunay property of $P$ and nature of 
$\cQ$}. More generally, let $\textsf{F}$ be a finite-dimensional linear  space of $\R$-valued 
functions  on $\Lambda$. Then a set $R \subset \Lambda$ is called perfect with respect to 
$\textsf{F}$ if the system of \emph{ linear} inhomogeneous equations $\{f(\rr)=c\; \vline \; \rr 
\in R\}$  on the coefficients of $f$ has a unique solution in $\textsf{F}$ for any $c \neq 0$. 

We have a natural bijection between perfect affine quadratic lattices $\Aff(\Lambda, Q)$ and 
triples $(\Lambda,P,\rho^2)$, where  $P$ is a perfect Delaunay polyhedron and $\rho^2 \ge 0$ is 
the squared radius of its Delaunay ellipsoid. Thus, there are only finitely many arithemtically inequivalent perfect Delaunay polyhedra in each dimension, up to rescaling. Furthermore, since for perfect Delaunay polyhedra arithmetic equivalence implies isometry, there are only finitely many nonisometric perfect Delaunay polytopes in each dimension, up to rescaling. 
For $n=0$ there is only one perfect Delaunay and it is
$\0$. For $n=1$ the polytope $[\0,\mathbf{1}]$ is perfect and Delaunay in $(\Z^1,x^2)$, and it 
is unique up to arithmetic equivalence.

\section{\protect \large Laminar Structure of Delaunay Cells}\label{sec:structure} 
 It is easy to see that a section of the vertex set of a Delaunay polytope by a 
 rational affine subspace is the vertex set 
of a Delaunay polytope in the induced sublattice. This observation suggests a 
recursive approach  to studying  Delaunay polytopes where each newly discovered  
Delaunay polytope is represented as a union of Delaunay polytopes of smaller 
dimensions lying in parallel subspaces. Indeed, for $n>1$  such a representation is 
never unique. It has been observed that dealing with a smaller numbers of big laminae  
is easier than studying a large number of small laminae. In other words, one is 
usually working with a representation in which the number of laminae is as small as 
possible. 

\begin{definition}
The lamina number $l(P)$ of a lattice polytope $P$ in a lattice $\Lambda$ is 
the minimal number of disjoint affine subspaces of $\Lambda \otimes \R$ whose 
intersections with $\vertex P$ form a partition of $\vertex P$ into proper subsets.  
\end{definition}   

The natural question is what laminar constructions lead to perfect Delaunay 
polytopes. In particular, is it possible to construct perfect Delaunay polytopes by 
using non-trivial (of dimension greater than 1) lower-dimensional perfect Delaunay 
polytopes  as some of the laminae? It turns out for $n=6-8$ such constructions are 
rather common, although not all of these polytopes have sections that are non-trivial 
perfect polytopes of smaller dimensions. The only perfect Delaunay $6$-polytope, 
Gosset's $G_6$, does not have non-trivial perfect sections. The only two known 
perfect $7$-polytopes, Gosset's $G_7$ and the 35-tope found by Erdahl and Rybnikov, 
each have a section isometric to $G_6$.  Our study showed that of the 27 known 
perfect Delaunay $8$-polytopes  17 have a section which is isometric to $G_6$, of 
which 10 have a section isometric to the 35-tope and one has a section isometric to 
$G_7$. The remaining 10 perfect 8-polytopes do not have non-trivial perfect sections.  

The lamina number $l(P)$ is closely related to the notion of lattice width. 
Denote by $\Lambda^{*}_{\cQ} \subset \aff \Lambda$ the dual of $\Lambda$ with respect to 
the  bilinear form  $\cQ$:  $\Lambda_{\cQ}^*$ consists of all vectors of $\aff 
\Lambda$ whose $\cQ$-products with vectors of $\Lambda$ are integer.
 If $B$ is a convex body in $\aff \Lambda$, then 
the width $w(B)$ of $B$ with respect to $(\Lambda,\cQ)$ is defined as the minimal 
natural number $w$ such that $B$ lies between hyperplanes $\cQ(\aaa^*,\x)=k$ and 
$\cQ(\aaa^*,\x)=k+w$, for some $\aaa^* \in \Lambda^*$.  It is widely believed (see e.g. 
Barvinok, 2002) that a body in $\aff \Lambda$ whose interior is empty of $\Lambda$-points 
cannot have lattice width exceeding $\rank \Lambda$.

\begin{proposition}
If $P$ is a Delaunay polytope, then $l(P)=w(P)+1$.
\end{proposition}

We do not know of any Delaunay polytopes whose lattice width exceeds 2. On the other hand, we have 

\begin{theorem}
If $P$ is a perfect (need not be Delaunay) polytope of dimension $n>1$, then $w(P)+1=l(P)>2$.
\end{theorem}
\begin{proof} Let $(\Lambda,\cQ)$ be the lattice in which $P$ is perfect.  Since the partition into laminae must be proper, $l(P)>1$. If $l(P)=2$, then there are affine sublattices $L_1$ and $L_2$  of codimension 1 such that $\vertex P =(\vertex P \cap L_1) \sqcup (\vertex P \cap L_1)$. There exists  an affine function $A$ on $\Lambda$, which is 1 on $L_1$ and  0 on $L_2$. Then the quadric $A\x\,(A\x-1)=0$ is circumscribed about $P$. If $Q$ is the quadratic function defined by $P$ uniquely up to scale, then $Q+A(A-1)$ must be proportional to $Q$, which means that $\Quad Q$ is of rank 1. Since $P$ is perfect, $\dim P=\rank \Lambda$. Since $P$ is a polytope,  $\rank \Lambda=
\rank\: \Quad Q=1$ and $\dim P=1$. 
\end{proof}

It turns out that all perfect Delaunay polytopes in dimensions $n=6-8$ have the lamina number $l$ equal to 3.
\begin{theorem}
Each perfect Delaunay polytope $P$ described in Sections 6--8 has $l(P)=3$  
\end{theorem}
\begin{proof} $\vertex G_6$ has a 3-laminae partition into a vertex, a 5-half-cube, and a 5-cross-polytope (Erdahl, Rybnikov 2002); another partition is into a 5-simplex, a 15-vertex polytope, and another 5 simplex  (see Erdahl , 1992). The partition of  $\vertex G_7$ 
into  the union of  the vertex sets of a 6-half-cube and two 6-cross-polytopes is given in 
Lemma   7.1.  The partition of $\Upsilon^7$  is given in 
Lemma   7.2  The partitions of the vertex sets of perfect 8-polytopes into layers follow from their coordinate representation given in Section 8.

\end{proof}
\subsection{Structure of Perfect Affine Lattices}\label{subsec:layer}
Recall that a pair $\Aff(\Lambda,Q)$, where $\Lambda$ is a lattice and $Q:\Lambda \rightarrow 
\R$ a quadratic function, is called perfect if the coefficients of $Q$ are uniquely determined 
from equations $X(\m)=\min\{Q(\z)\;\vline \; \z \in \Lambda\}$, where $\m$ runs over all 
minimal vectors of $Q$ and $X$ is an unknown quadratic function on $\Lambda$. In general, for a 
function $F$ defined on a lattice $\Lambda$ denote by $\cV(F)$ (the variety of $F$) the set of 
lattice points where $F$ is 0.  Perfection is a very natural notion as illustrated by the 
following theorem of Erdahl (1992). Recall that if $L$ is an affine sublattice of a lattice  
$\Gamma$, then $\overrightarrow{L}$ \textit{stands for the lattice}  $L-L$.
\begin{theorem}
$(\Lambda,Q)$ is perfect if and only if $\cV(Q)=\{\vv + \z \; \vline \; \vv \in \vertex P,\: \z 
\in \Gamma\}$, where $P$ is a perfect Delaunay polytope  in $\Aff(
\overrightarrow{\Lambda \cap \aff P},Q|_{\aff P})$ and $\Gamma$ is a sublattice of $\Lambda$ such that $\Lambda$ is the direct 
affine sum of $\Lambda \cap \aff P$ and  $\Gamma$ over $\Z$, i.e. 
\[\Lambda =\{(\x-\x') + \z \;\; \vline \;\; \x,\x' \in \Lambda \cap \aff P, \: \z \in \Gamma\}.\] 
\end{theorem} On the basis of this characterization Erdahl and Rybnikov  proved the following 
theorem (Erdahl, Ordine, Rybnikov, 2004).

\begin{theorem}\label{thm:step}Let $P$ be a perfect polytope in $Del(\Lambda,\cQ)$ and let $Q_P$ 
be its perfect quadratic function, i.e. $\vertex P=\cV(Q_P)$ and $\cQ=\Quad Q_P$. Suppose $D \in 
Del(\Lambda,\cQ)$ is another Delaunay cell of full dimension, which is not a 
$\Lambda$-translate of $P$. If $\e \notin \Lambda$, then there is a positive definite form 
$\cQ'$ on $\Lambda \oplus \Z \e$ and a perfect polytope  $P'$ in $Del(\Lambda \oplus \Z \e, 
\cQ')$ such that $P'\cap \aff \Lambda=P$ and $P' \cap \{\left(\aff \Lambda\right)+\e\}=D+\e$
\end{theorem}

The delicate part here is the case  where $P$ and $D$ have identical Delaunay radii. By using 
the following refinement of Erdahl's (1992) theorem it is possible to prove that under the 
assumptions of Theorem \ref{thm:step} a perfect polytope $P'$ that contains as sections $P$ and 
$D$ cannot be unbounded (Erdahl, Ordine, Rybnikov, 2004). On the basis of Theorem 
\ref{thm:step} we can make a few useful observations. 

\begin{itemize} \item If $P$ is an antisymmetric perfect polytope in $Del(\Lambda,\cQ)$, then 
there is a perfect polytope $P'$ in  $Del(\Lambda \oplus \Z\e,\cQ')$, for some form $\cQ$, with 
a section isometric to $P$. For example, for $P=G_6$  (Gosset's 6-polytope) there are two 
perfect 7-polytopes that have $G_6$ as a section, namely $G_7$ and the 35-tope.  $G_7$ can be  
obtained by taking $D=-P$ in Theorem \ref{thm:step}, while the 35-tope cannot be obtained by a 
direct application of Theorem \ref{thm:step}. \item We know only one example of a Delaunay 
tiling formed by translates of a \emph{centrally-symmetric} perfect Delaunay polytope: this is 
the tiling of $\Z^1$ by unit intervals. Incidently, we do not know of any Delaunay polytope, 
except for the $n$-cube, that tiles $\R^n$ by translation. We conjecture that for $n>1$ there 
are no such examples. If this is true, then Theorem \ref{thm:step} gives a construction for a 
new perfect Delaunay polytope in dimension $n+1$ from a perfect Delaunay polytope in dimension 
$n$. However, this construction is not uniquely defined, since  in Theorem \ref{thm:step} there 
may be  different choices of $D$.  \item When $P$ is a centrally symmetric perfect $n$-polytope 
in $Del(\Lambda,\cQ)$  and there is an antisymmetric $n$-polytope $D$  in  $Del(\Lambda,\cQ)$, 
it may happen that the center of the perfect polytope $P'$ coincides with the center of $P$. 
Then $P'$ has at least three $n$-dimensional layers, which are translates of $-D$, $P$, and $D$ 
respectively. The only 8-dimensional polytope from our list that has a section isometric to 
$G_7$ arises from  this construction. The role of $D$ is played by a Delaunay simplex of double 
volume in the Delaunay tiling defined by $G_7$  (see Erdahl, 1992 for a description). \item For 
many a perfect 8-polytope  the Delaunay tiling has a significant number of arithmetically 
inequivalent $8$-cells. This suggests that starting from $n=9$ the number of perfect Delaunay 
$n$-polytopes explodes.  (see \url{http://www.liga.ens.fr/~dutour} for the enumeration )  It is 
likely that $n=8$ is the highest dimension in which a complete classification is within reach. 
\end{itemize} \section{Symmetries of Perfect Delaunay Polytopes}

Recall that the group $O(\Z^n,\cQ)$ of  linear automorphisms of a \textit{quadratic} lattice 
$(\Z^n,\cQ)$ is defined as the 
full subgroup of $O(\R^n,\cQ)$ that maps  $\Z^n$ onto itself, in other words, the
 set-stabilizer of $\Z^n$ in  $O(\R^n,\cQ)$:
 \[O(\Z^n,\cQ)=\{\tau \in O(\R^n,\cQ) \;\; \vline \;\; \tau (\Z^n)=\Z^n\}.\] 
   The group $O(\Z^n,\cQ)$ can also 
be seen as the  subgroup of $GL_n(\Z)$ that consists of transformations preserving $\cQ$:
\[O(\Z^n,\cQ)=\{\tau \in GL(n,\Z) \;\; \vline \;\; \forall \z \in \Z^n: \cQ[\tau\z]=\cQ[\z]\}.\]
   Denote by $Iso(\R^n,\cQ)$ the  group of affine automorphisms of $\R^n$ which preserve  $\cQ$. If $D$ is a 
$\Z^n$-polytope, then $Iso(D,\cQ)$ denotes the group of all transformations 
from $Iso(\R^n,\cQ)$ that map $D$ onto itself.  Denote by $LatIso(D,\cQ)$  the group of all 
transformations from $O(\Z^n,\cQ)$ that map $D$ to itself. Clearly $LatIso(D,\cQ) \le Iso(D,\cQ)$. When  $\{\x-\y\;\vline\;\x,\y \in \vertex D \}=\Z^n$ the polytope $D$ is called generating. All 
known perfect Delaunay polytopes are generating. Obviously, for generating polytopes 
$Iso(D,\cQ)=LatIso(D,\cQ)$. 

  An important invariant of an extreme Delaunay polytope $D$ in $(\Z^n,\cQ)$ is the dimension of 
the subspace of quadratic forms in $n$ variables preserved by $O(\Z^n,\cQ)$. We 
denote this space by $\QuadInv[D]$. 

The metric geometry of a Delaunay polytope $D$ is reflected in the \emph{norm 
spectrum} of $D$, which is just the set of all possible value for $\cQ[\x-\y]$, where 
$\x, \y \in \vertex D$. We denote the norm spectrum by $\Spec(D)$.

We have classified the isometry groups of all known perfect Delaunay polytopes for $n 
\le 8$. The isometry groups of six- and seven-dimensional perfect polytopes are 
distinct. Among the isometry groups of the 27 8-polytopes there are 21 non-isomorphic. 
Polytopes in the following five groups have isomorphic groups: \#2 and \#5; \#3 and 
\#13; \#12 and \#21; \#14, \#19, and \#25; \#24 and \#27. The most interesting is the 
case of \#2 and \#5: both polytopes have 72 vertices in two orbits of size 56. Their 
group contains $\cS_8$ as a subroup of index 2.

\section{Perfect Affine Quadratic Lattices for $n<6$}\label{sec:list:6}
$\Aff(\mathbf{0},0)$ is a perfect affine lattice of  rank $0$. 
All perfect affine lattices of rank 1 are obviously equivalent, up to scaling, to 
$\Aff(\Z^1,(x-\frac{1}{2})^2)$. The inequality $(x-\frac{1}{2})^2 \le \frac{1}{4}$ 
describes the Delaunay ellipsoid for the Delaunay polytope $[0,1]$. The Delaunay 
tiling for $(\Z^1,\Quad (x-\frac{1}{2})^2)=(\Z^1,x^2)$ consists of points of $\Z^1$, 
which are 0-dimensional polytopes, and segments $[k,k+1]$, where $k \in \Z$, which 
are 1-dimensional polytopes of the tiling. The symmetry group of $[0,1]$ consists of 
2 elements and is generated by reflection about $\frac{1}{2}$. Surprisingly, there 
are no perfect affine modules of ranks 2, 3, 4, and 5. This was first proven by Erdahl 
in 1975 (see also Erdahl, 1992).

\section{Perfect affine quadratic lattices of rank 6}

The affine lattice $\Aff(\Z^6,\cE_6[\x-\cc])$, where $\cE_6$ is given by
\begin{equation*} 
\setlength{\extrarowheight}{1pt} \cE_6(\x)= \x^t \: \begin{array}{|c|c|c|c|c|c|}
\hline  4 & 3 & 3 & 3 & 3 & 5\\
\hline 3 & 4 & 3 & 3 & 3 & 5\\
\hline 3 & 3 & 4 & 3 & 3 & 5\\
\hline 3 & 3 & 3 & 4 & 3 & 5\\
\hline 3 & 3 & 3 & 3 & 4 & 5\\
\hline  5 & 5 & 5 & 5 & 5 & 8\\
 \hline
\end{array}
\: \x \quad \hbox{and} \quad 
\cc=\frac{1}{3}\:\begin{array}{|r|}\hline1\\ \hline1\\ \hline1\\ \hline1\\ \hline1\\ \hline-2\\ \hline
\end{array}\:,
\end{equation*}
turns out to be perfect, which was first observed by Erdahl (1975). Quadratic form 
$\cE_6$ is  of type $E_6$, that is $(\Z^6,\cE_6[\x])$ is equivalent, up to scaling, 
to $(E_6, \sum x_i^2)$, where $E_6$ is a well-known root lattice in $\R^6$ (see, e.g. 
Coxeter, 1991).  Lattice $(E_6, \sum x_i^2)$ has first been constructed by Korkin and 
Zolotarev (1873), to which they referred as the \emph{fifth perfect form in 6 
variables} and denoted it by $X$.

Inequality $\cE_6[\x-\cc] \le \frac{4}{3}$ defines the Delaunay ellipsoid for a Delaunay 
polytope in $(\Z^6,\cE_6)$. The set of 27 vertices of this polytope is given by the following 
table ($[-1,0^{4};-1]$ means the entry $0$ is repeated $4$ times and all permutations of the 
first 5 positions are taken -- the last entry is separated by semicolumn and is not permuted; 
the other records of the table are interpreted similarly). 
\begin{center} 
\setlength{\extrarowheight}{4pt} \begin{tabular}{|c|c|c|c|c|} \hline $x_6=-3$  & $x_6=-2$ &  
$x_6=-1$ &$x_6=0$ & $x_6=1$ \\ \hline  \hline  

$ \mathbf{[1^{5};-3]}\times 1$   & $\mathbf{[ 0,1^4;-2]}\times 5$  & 
$\mathbf{[1^{2},0^3;-1]}\times 10$ & $\mathbf{[0^{6}]} \times 1$& $\mathbf{[1,0^{4};-1]}\times 
5$
\\ \hline  
 &  &  & $\mathbf{[1,0^4;0]}\times 5$ & \\ \hline
\end{tabular}
\end{center}

This polytope is  known as Gosset's (1900) 6-dimensional semiregular polytope, which we denote 
by $G_6$. It is a two distance set and $\Spec(G_6)=\{2,4\}$. Typically, $G_6$ is described as a 
Delaunay polytope for $(E_6, \sum x_i^2)$, in which case it is commonly denoted by $2_{21}$ -- 
the notation going back to Coxeter (see his 1988 paper for history).  $G_6$ has two orbits of 
facets, regular simplices and regular cross-polytopes.  It does not have interior diagonals and 
all segments joining its vertices are either edges, or diagonals of its facets. The 1-skeleton 
of $G_6$ is a strongly regular graph known as Schlafli graph. Its isometry group 
$Iso(G_6,\cE_6)$, of order $51820$,  is the famous group of automorphisms of the 27 lines on a 
general cubic  surface. $Iso(G_6,\cE_6)$  is isomorphic to the semidirect product of a 
2-element group generated by a reflection and a reflection-free normal that consists of  
$2^{6}\times 3^{4} \times 5=25920$ elements; the latter group is simple and has a number of 
descriptions as a group of Lie type (see ATLAS, 1986, for more details). The isometry group of 
$G_6$ is transitive on its vertex set.  A remarkable property of $G_6$ is that the convex hull 
of any subset of $\vertex G_6$ is a Delaunay polytope for $(\Z^6,\cQ)$ for some positive definite 
form $\cQ$ (Erdahl and Rybnikov, 2002). For more details see Coxeter (1988) and Erdahl \& 
Rybnikov (2002).

\subsubsection{Laminar Structure of Gosset's $G_6$}
Let us denote by $J(n,s)$ the polytope formed by all $\{0,1\}$-vectors in $(\Z^n,\sum^{n}_{i=1} x_i)$ with the sum of the coordinates equal to $s$. It is known that for each $s$, such that   $0 \le s<n$, the polytope $J(n,s)$ is isometric to a  Delaunay
polytopes in $(A_{n-1}, \sum^{n-1}_{i=1} x_i)$, where $A_{n-1}$ is a root lattice of type $A$ of rank $n-1$.  $G_6$ can be represented as the union of 3 laminae that are isometric to $ J(6,1), \:J(6,2)$, and $J(6,1)$, where the regular 6- simplices ( $J(6,1)$) are parallel (see Erdahl, 1992). Another lamination of $\vertex G_6$ is into 3 layers that consist of 
a 0-simplex, a 5-halfcube, and a $5$-cross-polytope (see Erdahl, Rybnikov, 2002). These 
laminations correspond to the  subdiagrams of types  $A_5$ and $D_5$ of the Coxeter diagram $E_6$, which represents
 the isometry group of $G_6$ (see Humphreys, 1990).

It was long suspected that $\Aff(\Z^6,\cE_6[\x-\cc])$ is the only perfect affine lattice of 
rank 6 up to scaling. Finally, Dutour (2004), using his \emph{EXT-HYP7} program, 
created in 2002, proved that this is the case.

\section{Perfect Affine Quadratic Lattices of rank 7}\label{sec:list7} \subsection{Gosset 
Polytope in Lattice $E_7$} The affine lattice $\Aff(\Z^7,\cE_7[\x-\cc])$, where $\cE_7$ is 
given by

\begin{equation*} \setlength{\extrarowheight}{1pt} \cE_7(\x)=\x^t\; 
\begin{array}{|c|c|c|c|c|c|c|}
\hline4 & 3 & 3 & 3 & 3 & 5 & 4 \\
\hline3 & 4 & 3 & 3 & 3 & 5 & 4 \\
\hline3 & 3 & 4 & 3 & 3 & 5 & 4 \\
\hline3 & 3 & 3 & 4 & 3 & 5 & 4 \\
\hline3 & 3 & 3 & 3 & 4 & 5 & 4 \\
\hline5 & 5 & 5 & 5 & 5 & 8 & 6 \\
\hline4 & 4 & 4 & 4 & 4 & 6 & 6 \\
\hline \end{array} \;\x \quad \hbox{and} \quad 
\cc=\frac{1}{2}\;\begin{array}{|r|}\hline0\\\hline0\\\hline0\\\hline0\\\hline0\\\hline0\\\hline1\\\hline 
\end{array} \end{equation*} turns out to be perfect, which was first observed by Erdahl (1975). 
Quadratic form $\cE_7$ is  of type $E_7$, that is $(\Z^7,\cE_7)$ is equivalent, up to scaling, 
to $(E_7, \sum x_i^2)$, where $E_7$ is a well-known root lattice in $\R^7$ (see, e.g. Coxeter, 
1991). Lattice $(E_7, \sum x_i^2)$ has first been constructed by Korkine and Zolotareff (1873), 
to which they referred as the \emph{sixth perfect form in 7 variables} and denoted it by $Y$.

Inequality $\cE_7[\x-\cc] \le 3$ defines the Delaunay ellipsoid for a polytope in 
$Del(\Z^7,\cE_7)$, whose vertex set is given below. 

\begin{center} \setlength{\extrarowheight}{4pt} \begin{tabular}{|c|c|c|c|}\hline {$x_7=-1$} & 
{$x_7=0$} & {$x_7=1$}  & {$x_7=2$}\\ \hline \hline $\mathbf{[1^5;-2;-1]} \times 1$ & $\mathbf{[ 
0^{7}]}\times 1$ & $\mathbf{[ 0^{6};1]}\times 1$& $\mathbf{[-1^5;2;2]}  \times 1$\\ \hline
 & $\mathbf{[1,0^4;0;0]}\times 5$ &$\mathbf{[-1,0^4;0;1]}\times 5$ & \\
\cline{2-3}&$\mathbf{[-1,0^{4};1;0]}\times 5$  & $\mathbf{[1,0^{4};-1;1]}\times 5$ & \\ 
\cline{2-3}&$\mathbf{[1^{2},0^3;-1;0]}\times 10$& $\mathbf{[-1^{2},0^3;1;1]}\times 10$&\\ 
\cline{2-3}&$\mathbf{[0,1^4;-2;0]}\times 5$  & $\mathbf{[0,-1^4;2;1]}\times 5$ &\\ \cline{2-3} 
&$\mathbf{[1^{5};-3;0]}\times 1$ & $\mathbf{[-1^{5};3;1]}\times 1$&\\ \hline
 \end{tabular} \end{center}

This centrally-symmetric polytope has 56 vertices and is known as Gosset's (1900) 7-dimensional 
semiregular polytope, which we denote by $G_7$. It is a 3-distance set and 
$\Spec(G_7)=\{2,4,6\}$. Typically, $G_7$ is described as a Delaunay polytope for $(E_7, \sum 
x_i^2)$, in which case it is commonly denoted by $3_{21}$ after Coxeter (see Coxeter (1988) for 
history).  The 1-skeleton  of $G_7$ is known as Gosset graph, which is a strongly-regular 
graph. $G_7$ has 28 interior diagonals passing through its center; in fact, Patrick du Val 
discovered that  $G_7$ can be thought of as the convex hull of seven congruent 3-dimensional 
cubes in $\R^7$ with common center (see Coxeter, 1988). The isometry group of $G_7$ is 
transitive on vertices and is isomorphic to the semidirect product of a 2-element group 
generated by a reflection and a reflection-free normal subgroup of $4 \times 9!= 1451520$ 
elements; the latter group is simple and is isomorphic, among other groups of Lie type, to 
$O_7(2)$ (see ATLAS, 1986). 

\begin{lemma} $l(G_7)=3$ and  $\vertex(G_7)$ is the union of the vertex sets of two 
6-cross-polytopes and a 6-half-cube. \end{lemma} \begin{proof} The fact $l(G_7)=3$ follows from 
the representation of $G_7$ as the following subset of the 7-cube $[-1,+1]^7$. Consider all 
cyclic permutations of 8 vectors $(\pm1,\pm1;0;\pm1;0,0,0)$. This defines a 56-element subset 
$V$ of $[-1,+1]^7$.  P. Du Val (attributed to du Val by Coxeter, 1988) had shown that the convex 
hull of these 56 points is the Gosset 7-polytope in the lattice $\aff_{\Z} V=\Z^7$ with respect 
to the usual metric $\sum x_i^2$. Note that  $\conv V \notin Del(\Z^7,\sum x_i^2)$, since  
$\conv V$  obviously contains the origin. 

In each of the coordinate directions $\conv V$ has three layers defined by inequalities 
$x_i=-1,0,+1$  and thus $l(\conv V)=l(G_7)=3$. It is easy to see that the sections $x_i=-1$ and 
$x_i=+1$ are 6-cross-polytopes and the section $x_i=0$ is a 6-half-cube. \end{proof} 

\subsubsection{Laminar Structure of Gosset's $G_7$} $G_7$ can be represented as the union of 4 
laminae that are isometric to $ J(7,1), \:J(7,2), \:J(7,2),$ and $J(7,1)$ (obviously, $J(7,1)$ 
is regular simplex).  Another lamination of $\vertex G_7$ is into 3 layers that consist of 
vertex sets of a 6-cross-polytope, 6-halfcube, and another 6-cross-polytope.  Yet another 
lamination is given in the above table, where the layers are a 0-simplex, a copy of  $G_6$, 
another copy of $G_6$, and another 0-simplex. All these laminations correspond to the 
\emph{unique} subdiagrams of types  $A_6,\: D_6$, and  $E_6$ of the Coxeter diagram $E_7$, 
which represents the isometry group of $G_7$ (see Humphreys, 1990). 

\subsection{The 35-tope} The only known perfect 
affine lattice of rank 7, that is not equivalent to $\Aff(\Z^7,\cE_7[\x-\cc])$, was constructed 
by Erdahl and Rybnikov in 2000 (see  Erdahl and Rybnikov (2002) and Rybnikov (2001)). It is 
$\Aff(\Z^7,\cER_7[\x-\cc])$, where $\cER_7$ is given by

\begin{equation*} \setlength{\extrarowheight}{1pt} \cER_7(\x)=\x^t\; 
\begin{array}{|c|c|c|c|c|c|c|}
\hline 8 & 6 & 6 & 6 & 6 & 6 & 9 \\
\hline 6 & 8 & 6 & 6 & 6 & 6 & 9 \\
\hline 6 & 6 & 8 & 6 & 6 & 6 & 9 \\
\hline 6 & 6 & 6 & 8 & 6 & 6 & 9 \\
\hline 6 & 6 & 6 & 6 & 8 & 6 & 9 \\
\hline 6 & 6 & 6 & 6 & 6 & 8 & 9 \\
\hline 9 & 9 & 9 & 9 & 9 & 9 & 13 \\
\hline
\end{array}
\;\x \quad \hbox{and}  \quad 
\cc=\frac{1}{16}\;\begin{array}{|r|}\hline5\\\hline5\\\hline5\\\hline5\\\hline5\\\hline5\\\hline-14\\\hline
\end{array}
\end{equation*}

The lattice $(\Lambda, {\cER_7})$ has 12 shortest vectors  
and $ \det \Lambda_{\cER_7}=256 $. The order of $O(\Lambda,{\cER_7})$ 
is 2880 and the dimension of the space of invariant forms is 3. $(\Lambda, {\cER_7})$ is 
not perfect, but the lattice obtained from $(\Lambda, {\cER_7})$ by adding the centers 
of perfect ellipsoids is perfect with 70 shortest vectors. 

Inequality $\cER_7[\x-\cc] \le \frac{43}{16}$ defines the Delaunay ellipsoid for a perfect 
Delaunay  $\Upsilon^7$ in $Del(\Z^7,\cER_7)$, whose vertex set is given below.

\begin{center} 
\setlength{\extrarowheight}{4pt} \begin{tabular}{|c|c|c|c|c|} \hline $x_7=-4$  & 
$x_7=-3$ &  $x_7=-1$ &$x_7=0$ & $x_7=1$ \\ \hline  \hline $ \mathbf{[1^{6};-4]}\times 1$   & 
$\mathbf{[ 0,1^5;-3]}\times 6$  & $\mathbf{[1^{2},0^4;-1]}\times 15$ & $\mathbf{[0^{7}]} \times 
1$ & $\mathbf{[-1,0^{5};1]}\times 6$ \\ \hline  
 &  &  & $\mathbf{[1,0^5;0]}\times 6$ &\\ \hline
\end{tabular}
\end{center}

$\Spec(\Upsilon^7)=\{3,4,5,7,8,9\}$. Polytope $\Upsilon^7$ generalizes to an infinite series of 
perfect Delaunay polytopes $\Upsilon^n$ ($n \ge 7$) with $\frac{n(n+3)}{2}$ vertices (see 
Erdahl, Rybnikov, 2002): 
\begin{center} 
\setlength{\extrarowheight}{4pt} 
\begin{tabular}{|c|c|c|}
 \hline $x_n=-(n-3)$  & $x_n=-(n-4)$ &  $x_n=-1$\\ \hline  \hline 
$ \mathbf{[1^{n-1};-(n-3)]}\times 1$   & $\mathbf{[ 0,1^{n-2};-(n-4)]} \times (n-1)$ & 
$\mathbf{[1^{2},0^{n-3};-1]}\times \frac{(n-1)(n-2)}{2}$ \\ \hline & 
$\mathbf{[1,0^{n-2};0]}\times (n-1)$ & \\ \hline
 \end{tabular}
\end{center} \begin{center} 
\setlength{\extrarowheight}{4pt} \begin{tabular}{|c|c|}
 \hline $x_n=0$ & $x_n=1$ \\ \hline  \hline 
$\mathbf{[0^{n}]} \times 1$ & $\mathbf{[-1,0^{n-2};1]}\times (n-1)$ \\ \hline 
$\mathbf{[1,0^{n-2};0]}\times (n-1)$ &\\ \hline
 \end{tabular}
\end{center}

 Polytope $\Upsilon^7$ has lamina number $l(\Upsilon^7)=3$ and can be represented as the union of Gosset polytope $G_6$ and 
regular simplices of dimensions 2 and 4 lying in parallel subspaces of $\R^7$:

\begin{lemma}
The vertex set of $\Upsilon^7$ is, up to scaling, isometric to  the disjoint union of 
the vertex sets of a Gosset polytope $G_6$, a regular 5-simplex, and a 1-simplex.
\end{lemma}
\begin{proof}
Consider the following partition of  
\[
\vertex \Upsilon^7=S_1 \bigsqcup S_2 \bigsqcup S_3=\] 
\[\{[-1,0^{5};1]\times 6, [0,1^{5};-3]\times 6, [1^{2},0^4;-1]\times 15\}\: 
\bigsqcup\: \{[0^{5},1;0] \times 6 \}\:  \bigsqcup\: \{[0^7],[1^6;-4]\}.
\]
Let us show that the affine rank of the first subset is 6. The first subset $S_1$ can 
be represented as  $S_{11} \bigsqcup S_{12} \bigsqcup S_{13}$, where 
$\{S_{11}=[-1,0^{5};1]\times 6\}$, $\{S_{12}=[0,1^{5};-3]\times 6\}$, $\{ 
S_{13}=[1^{2},0^4;-1]\times 15\}$.

Notice that \[[-1,0,0,0,0,0,1]=[0,1,1,0,0,0,-1]+([0,0,0,1,1,0,-1]-[1,1,1,1,1,0,-3])\] 

Applying cyclic permutations of the first six characters to the above identity, we see that, 
each element of $S_{11}$ can be written as $\p_1+(\p_2-\p)$, where $\p_1, \p_2 \in S_3$ and $\p 
\in S_2$. Therefore, $\aff (S_{11} \cup S_{12} \cup S_{13})=\aff (S_{12} \cup S_{13}) $.
 Since both $S_{12}$ and $S_{13}$ lie on the hyperplane $2(x_1+x_2+x_3+x_4+x_5+x_6)+3x_7=1$,   $\dim \aff (S_1 \cup S_2 \cup S_3)=6$. 
 Notice that this argument does not work if we replace $\Upsilon^7$ with $\Upsilon^n$ for 
$n>7$. It turns out that for $n>7$  we have $\aff (S_{11} \cup S_{12} \cup S_{13}) 
\neq \aff (S_{12} \cup S_{13}) $.
 
It is clear that the affine subspaces generated by $S_2$ and $S_3$ are parallel to 
that generated by the first subset. Computing squared distances between the elements 
of the first subset  with metric $\cER_7$ shows that it is a two-distance set 
isometric to $G_6$ (dilated by a factor of 2).  The second and the third sets are 
obviously regular simplices.
\end{proof}

\section{Perfect Affine Quadratic Lattices of Rank 8}\label{sec:list8} We denote the perfect  
quadratic functions by $\cQ^8_i[\x-\cc]$, where $\cc \in \Q^n$, and the corresponding perfect  
Delaunay polytopes by $D^8_i$. For each  $\cQ^8_i[\x-\cc]$ we give 

\begin{itemize}
 \item  an integer Gram matrix, 
 \item  the center $\cc$ of the perfect  
ellipsoid, 
\item  the order of the group $O(\Z^8,\cQ^8_i)$ (the group's generators are 
available from the first author upon request), together with the size of the maximal symmetric 
subgroup,
 \item the number $s(\Z^8,\cQ^8_i)$ of shortest vectors,
 \item the 
dimension  of the subspace of $\Sym(8,\R)$ that consists of forms $\cQ$ such that 
$\cQ[T\z]=\cQ[\z]$ for every  $T \in O(\Z^8,\cQ^8_i)$, we denote this subspace by $\QuadInv 
\left[ O(\Z^8,\cQ^8_i) \right]$. 
 \end{itemize}
 
 For each  $D^8_i$ we give 
 \begin{itemize}
 \item the coordinates of the vertices,
 
\item $|Iso(D^8_i,\cQ^8_i)|=\textrm{Order of the isometry group of } D^8_i $,

\item Whether $D^8_i$ is Centrally-symmetric or Antisymmetric, 

\item  Maximal non-trivial perfect polytope of smaller dimension (i.e. $G_6$, $G_7$, or the 
$35$-tope) contained in $D^8_i$,  

\item Information on certain types of Delaunay polytopes contained in $D^8_i$. If $X$ is an 
arithmetic type of a lattice Delaunay polytope such as, e.g., $J(n,s)$, and $D^8_i$ contains a 
copy of $J(n,s)$, which is not a proper subpolytope of a $J(n',s') \subset D^8_i$, then we 
state that $J(n,s)$ is maximally included into  $D^8_i$.

\item  the norm spectrum  $\Spec(D^8_i)$. 

\item the lamina number $l(D^8_i)$ (always 3). 
\end{itemize} 

 The coordinatization of all polytopes is chosen so that the three laminae 
structure is transparent.  For some of the polytopes we give additional geometric information. 

\subsection{Delaunay Tilings of Lattices $A_n$ and $D_n$}
The geometric structure of 8-dimensional perfect Delaunay polytopes can be analyzed by relating the geometry of these polytopes
to the geometry of Delaunay tilings of lattices $A_n$ and $D_n$, which is explicit in our 8-dimensional data. In this section let 
$\e_1,\ldots,\e_n$ be the standard basis of $\Z^n \subset \E^n$, and let $I=\{\x \in \R^n \;\; \vline \;\;
 \forall i: 0\le \x \cdot \e_i \le 1\}$ denote the  standard unit cube.

$A_n$ can be defined as 
$(\Z^n, \sum_{i=1}^n x_i^2 + \sum_{1 \le i < j \le n} x_ix_j)$ or, in terms of the Euclidean space $\E^n$, as the lattice
based on a regular $n$-simplex. Lattice $D_n$ can be defined as 
 $(\Z^n, \sum_{i=1}^n x_i^2 + x_1x_3+\sum_{2 \le i < j \le n} x_ix_j)$ or, in terms of the Euclidean space $\E^n$, as the sublattice 
 of $\Z^n$ that consists of all points with even sum of the coordinates; another Euclidean construction of $D_n$ is obtained 
 by taking $\Z^n$ and adding to it the centers of all facets of the unit $n$-cubes with integral vertices -- this is an $n$-dimensional 
 generalization of what is known in crystallography as the face-centered cubic lattice, or \textsf{fcc.} Note that for $n=3$, $A_n$ and $D_n$
 coincide.
 
 Delaunay tilings of $A_n$ have been described by Barnes (1959) and Delaunay tilings of $D_n$ have been described by Ryshkov and Shushbaev (1981).
 Below we give a brief description of these tilings borrowed from Baranovski (1991). 
 
Let $\bd=\sum \e_i$. Consider the sections 
 of the standard unit cube by hyperplanes perpendicular to $\bd$ and passing through the points $\frac{q}{n}\bd$ 
 for $q=1,\ldots,n$. These hyperplanes induce a partition of $I$ into $n$ $n$-polytopes $P(q)$, where each $P(q)$ is squeezed between hyperplanes 
 $\sum x_i=q-1$ and $\sum x_i=q$. It can be shown (see Barnes (1959) or Baranovski (1991)) that with respect to quadratic form 
$\sum_{i=1}^n x_i^2 + \sum_{1 \le i < j \le n} x_ix_j$ these polytopes are Delaunay. Thus, any Delaunay $n$-polytope in $A_n$ 
is a translate of one these polytopes. The 1-skeletons of the faces of $P(q)$ that are defined by 
$P(q) \cap \{\x \in \R^n \;\; \vline \;\;\sum x_i=q-1\}$ and $P(q) \cap \{\x \in \R^n \;\; \vline \;\;\sum x_i=q\}$
are Johnson graphs $J(n,q-1)$ and $J(n,q)$. We also use $J(n,q-1)$ and  $J(n,q)$ to refer to the arithmetic classes of
these polytopes.  Note that $J(n,q)$ and $J(n,n-q)$ are isometric with respect to any quadratic form on $\Z^n$, since one of them can be obtained
from the other by a combination of a lattice translation and an inversion with respect to a lattice point; in the terminology of geometry
of numbers such polytopes are known as \emph{homologous.}
 
Let us consider  $D_n$  as the sublattice 
 of $\Z^n \subset \E^n$ that consists of all points with even sum of the coordinates. Then any Delaunay $n$-polytopes 
 is homologous to one of the following:

 \begin{enumerate}              
 \item a cross-polytope, with vertices in $D_n$, centered at $\x \in \Z^n$, where $\sum x_i\equiv1 \mod 2$, 
 \item the convex hull of points of $D_n$ that belong to the standard cube $I$,
 \item the convex hull of points of $D_n$ that belong to the shifted cube $I+\e_n$.
 \end{enumerate}

The polytopes in 2) and 3) are known as $n$-semicubes (or halfcubes). Note that for $n=3$ a semicube is a tetrahedron and for $n=4$ a semicube is
a cross-polytope; the latter fact explains why the Delaunay tiling of $D_4$ is formed by three homology classes of cross-polytopes.

\renewcommand\cV{\vert} 
\renewcommand\cF{\cQ} 
\renewcommand\lambda{l} 
\renewcommand\c{\mathbf{c}} 
\medskip

\centerline{\large \textbf{8-dimensional perfect Delaunay polytopes}} 
\bigskip
 Perfect affine quadratic lattice $\Aff(\Z^8,\cQ^8_{1}[x-\cc])$, where $\cQ^8_{1}$ is given by
\begin{equation*}
\setlength{\extrarowheight}{1pt} \cQ^8_{1}(\x)=\x^t\; \begin{array}{|c|c|c|c|c|c|c|c|}
\hline  4 & 1 & 1 & 1 & 1 & 1 & 1 & -5\\
\hline  1 & 4 & 1 & 1 & 1 & 1 & 1 & -5\\
\hline  1 & 1 & 4 & 1 & 1 & 1 & 1 & -5\\
\hline  1 & 1 & 1 & 4 & 1 & 1 & 1 & -5\\
\hline  1 & 1 & 1 & 1 & 4 & 1 & 1 & -5\\
\hline  1 & 1 & 1 & 1 & 1 & 4 & 1 & -5\\
\hline  1 & 1 & 1 & 1 & 1 & 1 & 4 & -5\\
\hline  -5 & -5 & -5 & -5 & -5 & -5 & -5 & 19\\
\hline\end{array}
\;\x \quad \hbox{and}  \quad
\cc=\frac{1}{10}\;\begin{array}{|r|}
\hline 7\\
\hline 7\\
\hline 7\\
\hline 7\\
\hline 7\\
\hline 7\\
\hline 7\\
\hline 10\\
\hline\end{array}\;.
\end{equation*}
\begin{itemize}
\item $|O(\Z^8,\cQ^8_{1})|=\textrm{20160}$; $S_{7}<O(\Z^8,\cQ^8_{1})$
\vspace{-2.5mm} \item $s(\Z^8,\cQ^8_{1})=\textrm{14}$
\vspace{-2.5mm} \item $\dim \QuadInv[\Z^8,\cQ^8_{1}]=\textrm{3}$
\end{itemize}
Inequality $\cQ_{1}^8[\x-\cc] \le \frac{43}{10}$ defines the Delaunay ellipsoid for a perfect polytope $D^8_{1}\in Del(\Z^8,\cQ^8_{1})$, whose vertex set ($|\vertex D^8_{1}|=\textrm{44}$) is given below.
\begin{center}
\setlength{\extrarowheight}{4pt}
\begin{tabular}{|c|c|c|}
\hline$x_{8}=0$ &$x_{8}=1$ &$x_{8}=2$\\
\hline\hline$\mathbf{[0^{8}]} \times 1$ &$\mathbf{[0^{2},1^{5};1]}\times 21$ &$\mathbf{[1^{7};2]} \times 1$\\
\cline{1-1}\cline{2-2}\cline{3-3}$\mathbf{[0^{6},1;0]}\times 7$ &$\mathbf{[0,1^{6};1]}\times 7$ &$\mathbf{[1^{6},2;2]}\times 7$\\
\hline
\end{tabular}
\end{center}
\begin{itemize}
\item $\Spec(D^8_{1})=\textrm{4, 6, 7, 9, 10, 12, 13, 15}$
\vspace{-2.5mm}\item $|Iso(D^8_{1},\cQ^8_{1})|=\textrm{10080}$\vspace{-2.5mm}\item $l(D^8_{1})=\textrm{3}$
\vspace{-2.5mm}\item Antisymmetric
\vspace{-2.5mm}\item Maximally contained subpolytopes: $H(2)$, $\frac12 H(4)$, $J(8,6)$
\end{itemize}
 
 Perfect affine quadratic lattice $\Aff(\Z^8,\cQ^8_{
2}[x-\cc])$, where $\cQ^8_{2}$ is given by
\begin{equation*}
\setlength{\extrarowheight}{1pt} \cQ^8_{2}(\x)=\x^t\; \begin{array}{|c|c|c|c|c|c|c|c|}
\hline  8 & 6 & 6 & 6 & 6 & 10 & 8 & 4\\
\hline  6 & 8 & 6 & 6 & 6 & 10 & 8 & 5\\
\hline  6 & 6 & 8 & 6 & 6 & 10 & 8 & 6\\
\hline  6 & 6 & 6 & 8 & 6 & 10 & 8 & 4\\
\hline  6 & 6 & 6 & 6 & 8 & 10 & 8 & 4\\
\hline  10 & 10 & 10 & 10 & 10 & 16 & 12 & 7\\
\hline  8 & 8 & 8 & 8 & 8 & 12 & 12 & 7\\
\hline  4 & 5 & 6 & 4 & 4 & 7 & 7 & 7\\
\hline\end{array}
\;\x \quad \hbox{and}  \quad
\cc=\frac{1}{2}\;\begin{array}{|r|}
\hline 0\\
\hline 0\\
\hline 0\\
\hline 0\\
\hline 0\\
\hline 0\\
\hline 1\\
\hline 0\\
\hline\end{array}\;.
\end{equation*}
\begin{itemize}
\item $|O(\Z^8,\cQ^8_{2})|=\textrm{80640}$; $S_{8}<O(\Z^8,\cQ^8_{2})$
\vspace{-2.5mm} \item $s(\Z^8,\cQ^8_{2})=\textrm{16}$
\vspace{-2.5mm} \item $\dim \QuadInv[\Z^8,\cQ^8_{2}]=\textrm{2}$
\end{itemize}
Inequality $\cQ_{2}^8[\x-\cc] \le 3$ defines the Delaunay ellipsoid for a perfect polytope $D^8_{2}\in Del(\Z^8,\cQ^8_{2})$, whose vertex set ($|\vertex D^8_{2}|=\textrm{72}$) is given below.
\begin{center}
\setlength{\extrarowheight}{4pt}
\begin{tabular}{|c|c|c|}
\hline$x_{8}=-1$ &$x_{8}=0$ &$x_{8}=1$\\
\hline\hline$\mathbf{[-1^{2};0;-1^{2};2^{2};-1]} \times 1$ &$\mathbf{[-1^{5};2^{2};0]} \times 1$ &$\mathbf{[0;-1^{2};0^{2};1;0;1]} \times 1$\\
\cline{1-1}\cline{2-2}\cline{3-3}$\mathbf{[-1;0^{2};-1^{2};1;2;-1]} \times 1$ &$\mathbf{[-1^{5};3;1;0]} \times 1$ &$\mathbf{[0^{2};-1;0,1;0^{2};1]}\times 2$\\
\cline{1-1}\cline{2-2}\cline{3-3}$\mathbf{[-1;0^{2};-1^{2};2;1;-1]} \times 1$ &$\mathbf{[-1^{4},0;2;1;0]}\times 5$ &$\mathbf{[0^{7};1]} \times 1$\\
\cline{1-1}\cline{2-2}\cline{3-3}$\mathbf{[-1;0;1;0^{3};1;-1]} \times 1$ &$\mathbf{[-1^{2},0^{3};1^{2};0]}\times 10$ &$\mathbf{[1;0;-1;0^{4};1]} \times 1$\\
\cline{1-1}\cline{2-2}\cline{3-3}$\mathbf{[0^{6};1;-1]} \times 1$ &$\mathbf{[-1,0^{4};0;1;0]}\times 5$ &$\mathbf{[1;0^{2};1^{2};-2;0;1]} \times 1$\\
\cline{1-1}\cline{2-2}\cline{3-3}$\mathbf{[0^{2};1;-1,0;0;1;-1]}\times 2$ &$\mathbf{[-1,0^{4};1;0^{2}]}\times 5$ &$\mathbf{[1;0^{2};1^{2};-1^{2};1]} \times 1$\\
\cline{1-1}\cline{2-2}\cline{3-3}$\mathbf{[0;1^{2};0^{2};-1;1;-1]} \times 1$ &$\mathbf{[0^{8}]} \times 1$ &$\mathbf{[1^{2};0;1^{2};-2;-1;1]} \times 1$\\
\cline{1-1}\cline{2-2}\cline{3-3} &$\mathbf{[0^{6};1;0]} \times 1$ &\\
\cline{2-2} &$\mathbf{[0^{4},1;-1;1;0]}\times 5$ &\\
\cline{2-2} &$\mathbf{[0^{4},1;0^{3}]}\times 5$ &\\
\cline{2-2} &$\mathbf{[0^{3},1^{2};-1;0^{2}]}\times 10$ &\\
\cline{2-2} &$\mathbf{[0,1^{4};-2;0^{2}]}\times 5$ &\\
\cline{2-2} &$\mathbf{[1^{5};-3;0^{2}]} \times 1$ &\\
\cline{2-2} &$\mathbf{[1^{5};-2;-1;0]} \times 1$ &\\
\hline
\end{tabular}
\end{center}
\begin{itemize}
\item $\Spec(D^8_{2})=\textrm{3, 4, 5, 7, 8, 9, 12}$
\vspace{-2.5mm}\item $|Iso(D^8_{2},\cQ^8_{2})|=\textrm{80640}$\vspace{-2.5mm}\item $l(D^8_{2})=\textrm{3}$
\vspace{-2.5mm}\item Centrally-symmetric
\vspace{-2.5mm}\item Maximally contained subpolytopes: $35-tope$, $H(3)$, $\frac12 H(6)$, $J(9,7)$
\end{itemize}

 Perfect affine quadratic lattice $\Aff(\Z^8,\cQ^8_{
3}[x-\cc])$, where $\cQ^8_{3}$ is given by
\begin{equation*}
\setlength{\extrarowheight}{1pt} \cQ^8_{3}(\x)=\x^t\; \begin{array}{|c|c|c|c|c|c|c|c|}
\hline  11 & 8 & 8 & 8 & -3 & 4 & -20 & 4\\
\hline  8 & 11 & 8 & 8 & -3 & 4 & -20 & 4\\
\hline  8 & 8 & 11 & 8 & -3 & 4 & -20 & 4\\
\hline  8 & 8 & 8 & 11 & -3 & 4 & -20 & 4\\
\hline  -3 & -3 & -3 & -3 & 3 & -1 & 6 & -1\\
\hline  4 & 4 & 4 & 4 & -1 & 4 & -10 & 1\\
\hline  -20 & -20 & -20 & -20 & 6 & -10 & 48 & -10\\
\hline  4 & 4 & 4 & 4 & -1 & 1 & -10 & 4\\
\hline\end{array}
\;\x \quad \hbox{and}  \quad
\cc=\frac{1}{92}\;\begin{array}{|r|}
\hline 30\\
\hline 30\\
\hline 30\\
\hline 30\\
\hline 36\\
\hline -42\\
\hline 5\\
\hline -42\\
\hline\end{array}\;.
\end{equation*}
\begin{itemize}
\item $|O(\Z^8,\cQ^8_{3})|=\textrm{96}$; $S_{4}<O(\Z^8,\cQ^8_{3})$
\vspace{-2.5mm} \item $s(\Z^8,\cQ^8_{3})=\textrm{2}$
\vspace{-2.5mm} \item $\dim \QuadInv[\Z^8,\cQ^8_{3}]=\textrm{12}$
\end{itemize}
Inequality $\cQ_{3}^8[\x-\cc] \le \frac{189}{46}$ defines the Delaunay ellipsoid for a perfect polytope $D^8_{3}\in Del(\Z^8,\cQ^8_{3})$, whose vertex set ($|\vertex D^8_{3}|=\textrm{47}$) is given below.
\begin{center}
\setlength{\extrarowheight}{4pt}
\begin{tabular}{|c|c|c|}
\hline$x_{5}=-1$ &$x_{5}=0$ &$x_{5}=1$\\
\hline\hline$\mathbf{[0^{4};-1;0^{3}]} \times 1$ &$\mathbf{[-1,0^{3};0;-1^{3}]}\times 4$ &$\mathbf{[0^{4};1;-2;-1^{2}]} \times 1$\\
\cline{1-1}\cline{2-2}\cline{3-3} &$\mathbf{[0^{5};-2;-1;-2]} \times 1$ &$\mathbf{[0^{4};1;-1^{2};-2]} \times 1$\\
\cline{2-2}\cline{3-3} &$\mathbf{[0^{5};-2;-1^{2}]} \times 1$ &$\mathbf{[0^{4};1;-1^{3}]} \times 1$\\
\cline{2-2}\cline{3-3} &$\mathbf{[0^{5};-1^{2};-2]} \times 1$ &$\mathbf{[0^{3},1;1;0^{3}]}\times 4$\\
\cline{2-2}\cline{3-3} &$\mathbf{[0^{8}]} \times 1$ &$\mathbf{[0^{2},1^{2};1;-1;0;-1]}\times 6$\\
\cline{2-2}\cline{3-3} &$\mathbf{[0^{7};1]} \times 1$ &$\mathbf{[0,1^{3};1;0;1;0]}\times 4$\\
\cline{2-2}\cline{3-3} &$\mathbf{[0^{5};1;0^{2}]} \times 1$ &$\mathbf{[1^{5};-1;1;-1]} \times 1$\\
\cline{2-2}\cline{3-3} &$\mathbf{[0^{3},1;0;-1;0^{2}]}\times 4$ &$\mathbf{[1^{5};-1;1;0]} \times 1$\\
\cline{2-2}\cline{3-3} &$\mathbf{[0^{3},1;0^{3};-1]}\times 4$ &$\mathbf{[1^{5};0;1;-1]} \times 1$\\
\cline{2-2}\cline{3-3} &$\mathbf{[0^{3},1;0^{4}]}\times 4$ &\\
\cline{2-2} &$\mathbf{[0,1^{3};0^{2};1;0]}\times 4$ &\\
\hline
\end{tabular}
\end{center}
\begin{itemize}
\item $\Spec(D^8_{3})=\textrm{3, 4, 5, 6, 7, 8, 9, 10, 11, 12, 13, 14, 15}$
\vspace{-2.5mm}\item $|Iso(D^8_{3},\cQ^8_{3})|=\textrm{48}$\vspace{-2.5mm}\item $l(D^8_{3})=\textrm{3}$
\vspace{-2.5mm}\item Antisymmetric
\vspace{-2.5mm}\item Maximally contained subpolytopes: $G_6$, $H(3)$, $\frac12 H(5)$, $J(7,5)$
\end{itemize}

\include{Delaunay04_V4_GroupAttack}
 Perfect affine quadratic lattice $\Aff(\Z^8,\cQ^8_{
5}[x-\cc])$, where $\cQ^8_{5}$ is given by
\begin{equation*}
\setlength{\extrarowheight}{1pt} \cQ^8_{5}(\x)=\x^t\; \begin{array}{|c|c|c|c|c|c|c|c|}
\hline  15 & 11 & 11 & 11 & 11 & 11 & 11 & 11\\
\hline  11 & 9 & 8 & 8 & 8 & 8 & 8 & 8\\
\hline  11 & 8 & 9 & 8 & 8 & 8 & 8 & 8\\
\hline  11 & 8 & 8 & 9 & 8 & 8 & 8 & 8\\
\hline  11 & 8 & 8 & 8 & 9 & 8 & 8 & 8\\
\hline  11 & 8 & 8 & 8 & 8 & 9 & 8 & 8\\
\hline  11 & 8 & 8 & 8 & 8 & 8 & 9 & 8\\
\hline  11 & 8 & 8 & 8 & 8 & 8 & 8 & 9\\
\hline\end{array}
\;\x \quad \hbox{and}  \quad
\cc=\frac{1}{8}\;\begin{array}{|r|}
\hline -3\\
\hline -3\\
\hline -3\\
\hline -3\\
\hline -3\\
\hline -3\\
\hline -3\\
\hline -3\\
\hline\end{array}\;.
\end{equation*}
\begin{itemize}
\item $|O(\Z^8,\cQ^8_{5})|=\textrm{161280}$; $S_{8}<O(\Z^8,\cQ^8_{5})$
\vspace{-2.5mm} \item $s(\Z^8,\cQ^8_{5})=\textrm{56}$
\vspace{-2.5mm} \item $\dim \QuadInv[\Z^8,\cQ^8_{5}]=\textrm{2}$
\end{itemize}
Inequality $\cQ_{5}^8[\x-\cc] \le \frac{15}{8}$ defines the Delaunay ellipsoid for a perfect polytope $D^8_{5}\in Del(\Z^8,\cQ^8_{5})$, whose vertex set ($|\vertex D^8_{5}|=\textrm{72}$) is given below.
\begin{center}
\setlength{\extrarowheight}{4pt}
\begin{tabular}{|c|c|c|}
\hline$\sum x_i=-4$ &$\sum x_i=-3$ &$\sum x_i=-2$\\
\hline\hline$\mathbf{[2;-1^{6},0]}\times 7$ &$\mathbf{[-1^{3},0^{5}]}\times 56$ &$\mathbf{[-3;0^{6},1]}\times 7$\\
\cline{1-1}\cline{2-2}\cline{3-3}$\mathbf{[3;-1^{7}]} \times 1$ & &$\mathbf{[-2;0^{7}]} \times 1$\\
\hline
\end{tabular}
\end{center}
\begin{itemize}
\item $\Spec(D^8_{5})=\textrm{2, 3, 4, 5, 6}$
\vspace{-2.5mm}\item $|Iso(D^8_{5},\cQ^8_{5})|=\textrm{80640}$\vspace{-2.5mm}\item $l(D^8_{5})=\textrm{3}$
\vspace{-2.5mm}\item Antisymmetric
\vspace{-2.5mm}\item Maximally contained subpolytopes: $H(3)$, $\frac12 H(5)$, $J(8,5)$
\end{itemize}

 Perfect affine quadratic lattice $\Aff(\Z^8,\cQ^8_{
6}[x-\cc])$, where $\cQ^8_{6}$ is given by
\begin{equation*}
\setlength{\extrarowheight}{1pt} \cQ^8_{6}(\x)=\x^t\; \begin{array}{|c|c|c|c|c|c|c|c|}
\hline  8 & 6 & 6 & 6 & 6 & 6 & 9 & 5\\
\hline  6 & 8 & 6 & 6 & 6 & 6 & 9 & 5\\
\hline  6 & 6 & 8 & 6 & 6 & 6 & 9 & 5\\
\hline  6 & 6 & 6 & 8 & 6 & 6 & 9 & 4\\
\hline  6 & 6 & 6 & 6 & 8 & 6 & 9 & 6\\
\hline  6 & 6 & 6 & 6 & 6 & 8 & 9 & 5\\
\hline  9 & 9 & 9 & 9 & 9 & 9 & 13 & 7\\
\hline  5 & 5 & 5 & 4 & 6 & 5 & 7 & 6\\
\hline\end{array}
\;\x \quad \hbox{and}  \quad
\cc=\frac{1}{8}\;\begin{array}{|r|}
\hline 3\\
\hline 3\\
\hline 3\\
\hline 2\\
\hline 4\\
\hline 3\\
\hline -8\\
\hline -2\\
\hline\end{array}\;.
\end{equation*}
\begin{itemize}
\item $|O(\Z^8,\cQ^8_{6})|=\textrm{768}$; $S_{4}<O(\Z^8,\cQ^8_{6})$
\vspace{-2.5mm} \item $s(\Z^8,\cQ^8_{6})=\textrm{2}$
\vspace{-2.5mm} \item $\dim \QuadInv[\Z^8,\cQ^8_{6}]=\textrm{5}$
\end{itemize}
Inequality $\cQ_{6}^8[\x-\cc] \le \frac{11}{4}$ defines the Delaunay ellipsoid for a perfect polytope $D^8_{6}\in Del(\Z^8,\cQ^8_{6})$, whose vertex set ($|\vertex D^8_{6}|=\textrm{54}$) is given below.
\begin{center}
\setlength{\extrarowheight}{4pt}
\begin{tabular}{|c|c|c|}
\hline$x_{8}=-1$ &$x_{8}=0$ &$x_{8}=1$\\
\hline\hline$\mathbf{[0^{3};-1;1;0;1;-1]} \times 1$ &$\mathbf{[-1,0^{5};1;0]}\times 6$ &$\mathbf{[0^{7};1]} \times 1$\\
\cline{1-1}\cline{2-2}\cline{3-3}$\mathbf{[0^{6};1;-1]} \times 1$ &$\mathbf{[0^{8}]} \times 1$ &$\mathbf{[0^{3};1;-1;0^{2};1]} \times 1$\\
\cline{1-1}\cline{2-2}\cline{3-3}$\mathbf{[0^{4};1;0^{2};-1]} \times 1$ &$\mathbf{[0^{5},1;0^{2}]}\times 6$ &$\mathbf{[0^{3};1;0^{2};-1;1]} \times 1$\\
\cline{1-1}\cline{2-2}\cline{3-3}$\mathbf{[0^{2},1;0;1^{2};-1^{2}]}\times 3$ &$\mathbf{[0^{4},1^{2};-1;0]}\times 15$ &\\
\cline{1-1}\cline{2-2}$\mathbf{[0,1^{2};0;1;0;-1^{2}]}\times 3$ &$\mathbf{[0,1^{5};-3;0]}\times 6$ &\\
\cline{1-1}\cline{2-2}$\mathbf{[0,1^{2};0;1^{2};-2;-1]}\times 3$ &$\mathbf{[1^{6};-4;0]} \times 1$ &\\
\cline{1-1}\cline{2-2}$\mathbf{[1^{3};0;1;0;-2;-1]} \times 1$ & &\\
\cline{1-1}$\mathbf{[1^{3};0;2;1;-3;-1]} \times 1$ & &\\
\cline{1-1}$\mathbf{[1^{6};-3;-1]} \times 1$ & &\\
\cline{1-1}$\mathbf{[1^{4};2;1;-4;-1]} \times 1$ & &\\
\hline
\end{tabular}
\end{center}
\begin{itemize}
\item $\Spec(D^8_{6})=\textrm{2, 3, 4, 5, 6, 7, 8, 9, 10}$
\vspace{-2.5mm}\item $|Iso(D^8_{6},\cQ^8_{6})|=\textrm{384}$\vspace{-2.5mm}\item $l(D^8_{6})=\textrm{3}$
\vspace{-2.5mm}\item Antisymmetric
\vspace{-2.5mm}\item Maximally contained subpolytopes: $35-tope$, $H(3)$, $\frac12 H(5)$, $J(7,5)$
\end{itemize}

 Perfect affine quadratic lattice $\Aff(\Z^8,\cQ^8_{
7}[x-\cc])$, where $\cQ^8_{7}$ is given by
\begin{equation*}
\setlength{\extrarowheight}{1pt} \cQ^8_{7}(\x)=\x^t\; \begin{array}{|c|c|c|c|c|c|c|c|}
\hline  12 & 9 & 9 & 9 & 9 & -15 & 9 & 8\\
\hline  9 & 12 & 9 & 9 & 9 & -15 & 9 & 11\\
\hline  9 & 9 & 12 & 9 & 9 & -15 & 10 & 9\\
\hline  9 & 9 & 9 & 12 & 9 & -15 & 10 & 9\\
\hline  9 & 9 & 9 & 9 & 12 & -15 & 10 & 9\\
\hline  -15 & -15 & -15 & -15 & -15 & 24 & -15 & -14\\
\hline  9 & 9 & 10 & 10 & 10 & -15 & 12 & 8\\
\hline  8 & 11 & 9 & 9 & 9 & -14 & 8 & 13\\
\hline\end{array}
\;\x \quad \hbox{and}  \quad
\cc=\frac{1}{66}\;\begin{array}{|r|}
\hline 22\\
\hline 4\\
\hline 13\\
\hline 13\\
\hline 13\\
\hline 32\\
\hline 9\\
\hline 18\\
\hline\end{array}\;.
\end{equation*}
\begin{itemize}
\item $|O(\Z^8,\cQ^8_{7})|=\textrm{72}$; $S_{3}<O(\Z^8,\cQ^8_{7})$
\vspace{-2.5mm} \item $s(\Z^8,\cQ^8_{7})=\textrm{2}$
\vspace{-2.5mm} \item $\dim \QuadInv[\Z^8,\cQ^8_{7}]=\textrm{12}$
\end{itemize}
Inequality $\cQ_{7}^8[\x-\cc] \le \frac{91}{22}$ defines the Delaunay ellipsoid for a perfect polytope $D^8_{7}\in Del(\Z^8,\cQ^8_{7})$, whose vertex set ($|\vertex D^8_{7}|=\textrm{46}$) is given below.
\begin{center}
\setlength{\extrarowheight}{4pt}
\begin{tabular}{|c|c|c|}
\hline$x_{8}=-1$ &$x_{8}=0$ &$x_{8}=1$\\
\hline\hline$\mathbf{[0;1^{5};-1^{2}]} \times 1$ &$\mathbf{[-1,0^{4};-1;0^{2}]}\times 5$ &$\mathbf{[0;-1;-1^{2},0;-1;1^{2}]}\times 3$\\
\cline{1-1}\cline{2-2}\cline{3-3} &$\mathbf{[0^{5};-1^{2};0]} \times 1$ &$\mathbf{[0;-1;-1,0^{2};-1;0;1]}\times 3$\\
\cline{2-2}\cline{3-3} &$\mathbf{[0^{8}]} \times 1$ &$\mathbf{[0;-1;0^{4};1^{2}]} \times 1$\\
\cline{2-2}\cline{3-3} &$\mathbf{[0^{6};1;0]} \times 1$ &$\mathbf{[0^{7};1]} \times 1$\\
\cline{2-2}\cline{3-3} &$\mathbf{[0^{4},1;0^{3}]}\times 5$ &$\mathbf{[0^{5};1^{3}]} \times 1$\\
\cline{2-2}\cline{3-3} &$\mathbf{[0^{3},1^{2};1;0^{2}]}\times 10$ &$\mathbf{[1;-1;0^{5};1]} \times 1$\\
\cline{2-2}\cline{3-3} &$\mathbf{[0,1^{4};2;0^{2}]}\times 5$ &$\mathbf{[1;-1;0^{3};1^{3}]} \times 1$\\
\cline{2-2}\cline{3-3} &$\mathbf{[1^{5};2;-1;0]} \times 1$ &$\mathbf{[1;0^{4};1;0;1]} \times 1$\\
\cline{2-2}\cline{3-3} &$\mathbf{[1^{5};3;0^{2}]} \times 1$ &\\
\cline{2-2} &$\mathbf{[0^{2};1^{4};-1;0]} \times 1$ &\\
\cline{2-2} &$\mathbf{[0;1;0^{3};1^{2};0]} \times 1$ &\\
\cline{2-2} &$\mathbf{[1;0^{4};1^{2};0]} \times 1$ &\\
\hline
\end{tabular}
\end{center}
\begin{itemize}
\item $\Spec(D^8_{7})=\textrm{3, 4, 5, 6, 7, 8, 9, 10, 11, 12, 13, 15}$
\vspace{-2.5mm}\item $|Iso(D^8_{7},\cQ^8_{7})|=\textrm{36}$\vspace{-2.5mm}\item $l(D^8_{7})=\textrm{3}$
\vspace{-2.5mm}\item Antisymmetric
\vspace{-2.5mm}\item Maximally contained subpolytopes: $G_6$, $H(3)$, $\frac12 H(5)$, $J(7,5)$
\end{itemize}

 Perfect affine quadratic lattice $\Aff(\Z^8,\cQ^8_{
8}[x-\cc])$, where $\cQ^8_{8}$ is given by
\begin{equation*}
\setlength{\extrarowheight}{1pt} \cQ^8_{8}(\x)=\x^t\; \begin{array}{|c|c|c|c|c|c|c|c|}
\hline  8 & 6 & 6 & 6 & 6 & 6 & 9 & 5\\
\hline  6 & 8 & 6 & 6 & 6 & 6 & 9 & 6\\
\hline  6 & 6 & 8 & 6 & 6 & 6 & 9 & 6\\
\hline  6 & 6 & 6 & 8 & 6 & 6 & 9 & 5\\
\hline  6 & 6 & 6 & 6 & 8 & 6 & 9 & 5\\
\hline  6 & 6 & 6 & 6 & 6 & 8 & 9 & 7\\
\hline  9 & 9 & 9 & 9 & 9 & 9 & 13 & 8\\
\hline  5 & 6 & 6 & 5 & 5 & 7 & 8 & 8\\
\hline\end{array}
\;\x \quad \hbox{and}  \quad
\cc=\frac{1}{20}\;\begin{array}{|r|}
\hline 7\\
\hline 4\\
\hline 4\\
\hline 7\\
\hline 7\\
\hline 1\\
\hline -16\\
\hline 6\\
\hline\end{array}\;.
\end{equation*}
\begin{itemize}
\item $|O(\Z^8,\cQ^8_{8})|=\textrm{384}$; $S_{4}<O(\Z^8,\cQ^8_{8})$
\vspace{-2.5mm} \item $s(\Z^8,\cQ^8_{8})=\textrm{2}$
\vspace{-2.5mm} \item $\dim \QuadInv[\Z^8,\cQ^8_{8}]=\textrm{7}$
\end{itemize}
Inequality $\cQ_{8}^8[\x-\cc] \le \frac{14}{5}$ defines the Delaunay ellipsoid for a perfect polytope $D^8_{8}\in Del(\Z^8,\cQ^8_{8})$, whose vertex set ($|\vertex D^8_{8}|=\textrm{52}$) is given below.
\begin{center}
\setlength{\extrarowheight}{4pt}
\begin{tabular}{|c|c|c|}
\hline$x_{8}=-1$ &$x_{8}=0$ &$x_{8}=1$\\
\hline\hline$\mathbf{[0;1^{2};0^{2};1;-1^{2}]} \times 1$ &$\mathbf{[-1,0^{5};1;0]}\times 6$ &$\mathbf{[0;-1^{2};0^{2};-1;2;1]} \times 1$\\
\cline{1-1}\cline{2-2}\cline{3-3} &$\mathbf{[0^{8}]} \times 1$ &$\mathbf{[0;-1;0^{3};-1;1^{2}]} \times 1$\\
\cline{2-2}\cline{3-3} &$\mathbf{[0^{5},1;0^{2}]}\times 6$ &$\mathbf{[0^{2};-1;0^{2};-1;1^{2}]} \times 1$\\
\cline{2-2}\cline{3-3} &$\mathbf{[0^{4},1^{2};-1;0]}\times 15$ &$\mathbf{[0^{7};1]} \times 1$\\
\cline{2-2}\cline{3-3} &$\mathbf{[0,1^{5};-3;0]}\times 6$ &$\mathbf{[0^{3};0,1;-1;0;1]}\times 2$\\
\cline{2-2}\cline{3-3} &$\mathbf{[1^{6};-4;0]} \times 1$ &$\mathbf{[0^{3};0,1;0;-1;1]}\times 2$\\
\cline{2-2}\cline{3-3} & &$\mathbf{[0^{3};1^{2};-1^{2};1]} \times 1$\\
\cline{3-3} & &$\mathbf{[1;0^{4};-1;0;1]} \times 1$\\
\cline{3-3} & &$\mathbf{[1;0^{5};-1;1]} \times 1$\\
\cline{3-3} & &$\mathbf{[1;0^{2};0,1;-1^{2};1]}\times 2$\\
\cline{3-3} & &$\mathbf{[1;0^{2};1^{2};0;-2;1]} \times 1$\\
\cline{3-3} & &$\mathbf{[1;0;1^{3};0;-3;1]} \times 1$\\
\cline{3-3} & &$\mathbf{[1^{2};0;1^{2};0;-3;1]} \times 1$\\
\hline
\end{tabular}
\end{center}
\begin{itemize}
\item $\Spec(D^8_{8})=\textrm{2, 3, 4, 5, 6, 7, 8, 9, 10}$
\vspace{-2.5mm}\item $|Iso(D^8_{8},\cQ^8_{8})|=\textrm{192}$\vspace{-2.5mm}\item $l(D^8_{8})=\textrm{3}$
\vspace{-2.5mm}\item Antisymmetric
\vspace{-2.5mm}\item Maximally contained subpolytopes: $35-tope$, $H(3)$, $\frac12 H(5)$, $J(7,5)$
\end{itemize}

 Perfect affine quadratic lattice $\Aff(\Z^8,\cQ^8_{
9}[x-\cc])$, where $\cQ^8_{9}$ is given by
\begin{equation*}
\setlength{\extrarowheight}{1pt} \cQ^8_{9}(\x)=\x^t\; \begin{array}{|c|c|c|c|c|c|c|c|}
\hline  6 & 4 & 4 & 4 & 4 & 4 & -5 & 6\\
\hline  4 & 6 & 4 & 4 & 4 & 4 & -5 & 6\\
\hline  4 & 4 & 6 & 4 & 4 & 4 & -5 & 6\\
\hline  4 & 4 & 4 & 6 & 4 & 4 & -5 & 6\\
\hline  4 & 4 & 4 & 4 & 6 & 4 & -5 & 6\\
\hline  4 & 4 & 4 & 4 & 4 & 6 & -5 & 6\\
\hline  -5 & -5 & -5 & -5 & -5 & -5 & 9 & -6\\
\hline  6 & 6 & 6 & 6 & 6 & 6 & -6 & 9\\
\hline\end{array}
\;\x \quad \hbox{and}  \quad
\cc=\frac{1}{8}\;\begin{array}{|r|}
\hline 3\\
\hline 3\\
\hline 3\\
\hline 3\\
\hline 3\\
\hline 3\\
\hline 6\\
\hline -4\\
\hline\end{array}\;.
\end{equation*}
\begin{itemize}
\item $|O(\Z^8,\cQ^8_{9})|=\textrm{2880}$; $S_{6}<O(\Z^8,\cQ^8_{9})$
\vspace{-2.5mm} \item $s(\Z^8,\cQ^8_{9})=\textrm{12}$
\vspace{-2.5mm} \item $\dim \QuadInv[\Z^8,\cQ^8_{9}]=\textrm{5}$
\end{itemize}
Inequality $\cQ_{9}^8[\x-\cc] \le \frac{27}{8}$ defines the Delaunay ellipsoid for a perfect polytope $D^8_{9}\in Del(\Z^8,\cQ^8_{9})$, whose vertex set ($|\vertex D^8_{9}|=\textrm{58}$) is given below.
\begin{center}
\setlength{\extrarowheight}{4pt}
\begin{tabular}{|c|c|c|}
\hline$x_{7}=0$ &$x_{7}=1$ &$x_{7}=2$\\
\hline\hline$\mathbf{[-1,0^{5};0;1]}\times 6$ &$\mathbf{[0^{6};1^{2}]} \times 1$ &$\mathbf{[1^{6};2;-2]} \times 1$\\
\cline{1-1}\cline{2-2}\cline{3-3}$\mathbf{[0^{8}]} \times 1$ &$\mathbf{[0^{4},1^{2};1;0]}\times 15$ &\\
\cline{1-1}\cline{2-2}$\mathbf{[0^{7};1]} \times 1$ &$\mathbf{[0^{3},1^{3};1;-1]}\times 20$ &\\
\cline{1-1}\cline{2-2}$\mathbf{[0^{5},1;0^{2}]}\times 6$ &$\mathbf{[0,1^{5};1;-2]}\times 6$ &\\
\cline{1-1}\cline{2-2} &$\mathbf{[1^{7};-3]} \times 1$ &\\
\hline
\end{tabular}
\end{center}
\begin{itemize}
\item $\Spec(D^8_{9})=\textrm{3, 4, 5, 6, 7, 8, 9, 10, 11, 12}$
\vspace{-2.5mm}\item $|Iso(D^8_{9},\cQ^8_{9})|=\textrm{1440}$\vspace{-2.5mm}\item $l(D^8_{9})=\textrm{3}$
\vspace{-2.5mm}\item Antisymmetric
\vspace{-2.5mm}\item Maximally contained subpolytopes: $H(3)$, $\frac12 H(5)$, $J(7,5)$
\end{itemize}

 Perfect affine quadratic lattice $\Aff(\Z^8,\cQ^8_{
10}[x-\cc])$, where $\cQ^8_{10}$ is given by
\begin{equation*}
\setlength{\extrarowheight}{1pt} \cQ^8_{10}(\x)=\x^t\; \begin{array}{|c|c|c|c|c|c|c|c|}
\hline  8 & 6 & 6 & 6 & 6 & 6 & 9 & 4\\
\hline  6 & 8 & 6 & 6 & 6 & 6 & 9 & 4\\
\hline  6 & 6 & 8 & 6 & 6 & 6 & 9 & 3\\
\hline  6 & 6 & 6 & 8 & 6 & 6 & 9 & 3\\
\hline  6 & 6 & 6 & 6 & 8 & 6 & 9 & 2\\
\hline  6 & 6 & 6 & 6 & 6 & 8 & 9 & 3\\
\hline  9 & 9 & 9 & 9 & 9 & 9 & 13 & 4\\
\hline  4 & 4 & 3 & 3 & 2 & 3 & 4 & 5\\
\hline\end{array}
\;\x \quad \hbox{and}  \quad
\cc=\frac{1}{26}\;\begin{array}{|r|}
\hline 6\\
\hline 6\\
\hline 7\\
\hline 7\\
\hline 8\\
\hline 7\\
\hline -18\\
\hline 2\\
\hline\end{array}\;.
\end{equation*}
\begin{itemize}
\item $|O(\Z^8,\cQ^8_{10})|=\textrm{576}$; $S_{3}<O(\Z^8,\cQ^8_{10})$
\vspace{-2.5mm} \item $s(\Z^8,\cQ^8_{10})=\textrm{4}$
\vspace{-2.5mm} \item $\dim \QuadInv[\Z^8,\cQ^8_{10}]=\textrm{4}$
\end{itemize}
Inequality $\cQ_{10}^8[\x-\cc] \le \frac{35}{13}$ defines the Delaunay ellipsoid for a perfect polytope $D^8_{10}\in Del(\Z^8,\cQ^8_{10})$, whose vertex set ($|\vertex D^8_{10}|=\textrm{55}$) is given below.
\begin{center}
\setlength{\extrarowheight}{4pt}
\begin{tabular}{|c|c|c|}
\hline$x_{8}=-1$ &$x_{8}=0$ &$x_{8}=1$\\
\hline\hline$\mathbf{[1^{2};0;1;0;1;-2;-1]} \times 1$ &$\mathbf{[-1,0^{5};1;0]}\times 6$ &$\mathbf{[-1^{4};0^{2};3;1]} \times 1$\\
\cline{1-1}\cline{2-2}\cline{3-3}$\mathbf{[1^{3};0^{2};1;-2;-1]} \times 1$ &$\mathbf{[0^{8}]} \times 1$ &$\mathbf{[-1^{3};0^{2};-1;3;1]} \times 1$\\
\cline{1-1}\cline{2-2}\cline{3-3}$\mathbf{[1^{4};0^{2};-2;-1]} \times 1$ &$\mathbf{[0^{5},1;0^{2}]}\times 6$ &$\mathbf{[-1^{3};0^{3};2;1]} \times 1$\\
\cline{1-1}\cline{2-2}\cline{3-3}$\mathbf{[1^{4};0;1;-3;-1]} \times 1$ &$\mathbf{[0^{4},1^{2};-1;0]}\times 15$ &$\mathbf{[-1^{2};0;-1;0;-1;3;1]} \times 1$\\
\cline{1-1}\cline{2-2}\cline{3-3}$\mathbf{[1,2;1^{4};-4;-1]}\times 2$ &$\mathbf{[0,1^{5};-3;0]}\times 6$ &$\mathbf{[-1^{2};0;-1;0^{2};2;1]} \times 1$\\
\cline{1-1}\cline{2-2}\cline{3-3}$\mathbf{[2^{2};1^{2};0;1;-4;-1]} \times 1$ &$\mathbf{[1^{6};-4;0]} \times 1$ &$\mathbf{[-1^{2};0^{3};-1;2;1]} \times 1$\\
\cline{1-1}\cline{2-2}\cline{3-3}$\mathbf{[2^{2};1^{4};-5;-1]} \times 1$ & &$\mathbf{[-1^{2};0^{2};1;0;1^{2}]} \times 1$\\
\cline{1-1}\cline{3-3} & &$\mathbf{[-1,0;0^{4};1^{2}]}\times 2$\\
\cline{3-3} & &$\mathbf{[-1,0;0^{2};1;0^{2};1]}\times 2$\\
\cline{3-3} & &$\mathbf{[0^{7};1]} \times 1$\\
\hline
\end{tabular}
\end{center}
\begin{itemize}
\item $\Spec(D^8_{10})=\textrm{2, 3, 4, 5, 6, 7, 8, 9}$
\vspace{-2.5mm}\item $|Iso(D^8_{10},\cQ^8_{10})|=\textrm{288}$\vspace{-2.5mm}\item $l(D^8_{10})=\textrm{3}$
\vspace{-2.5mm}\item Antisymmetric
\vspace{-2.5mm}\item Maximally contained subpolytopes: $35-tope$, $H(3)$, $\frac12 H(5)$, $J(7,5)$
\end{itemize}

 Perfect affine quadratic lattice $\Aff(\Z^8,\cQ^8_{
11}[x-\cc])$, where $\cQ^8_{11}$ is given by
\begin{equation*}
\setlength{\extrarowheight}{1pt} \cQ^8_{11}(\x)=\x^t\; \begin{array}{|c|c|c|c|c|c|c|c|}
\hline  24 & 18 & 18 & 18 & 18 & 18 & 27 & 15\\
\hline  18 & 24 & 18 & 18 & 18 & 18 & 27 & 15\\
\hline  18 & 18 & 24 & 18 & 18 & 18 & 27 & 15\\
\hline  18 & 18 & 18 & 24 & 18 & 18 & 27 & 15\\
\hline  18 & 18 & 18 & 18 & 24 & 18 & 27 & 15\\
\hline  18 & 18 & 18 & 18 & 18 & 24 & 27 & 11\\
\hline  27 & 27 & 27 & 27 & 27 & 27 & 39 & 20\\
\hline  15 & 15 & 15 & 15 & 15 & 11 & 20 & 17\\
\hline\end{array}
\;\x \quad \hbox{and}  \quad
\cc=\frac{1}{90}\;\begin{array}{|r|}
\hline 36\\
\hline 36\\
\hline 36\\
\hline 36\\
\hline 36\\
\hline 22\\
\hline -91\\
\hline -21\\
\hline\end{array}\;.
\end{equation*}
\begin{itemize}
\item $|O(\Z^8,\cQ^8_{11})|=\textrm{480}$; $S_{5}<O(\Z^8,\cQ^8_{11})$
\vspace{-2.5mm} \item $s(\Z^8,\cQ^8_{11})=\textrm{4}$
\vspace{-2.5mm} \item $\dim \QuadInv[\Z^8,\cQ^8_{11}]=\textrm{8}$
\end{itemize}
Inequality $\cQ_{11}^8[\x-\cc] \le \frac{124}{15}$ defines the Delaunay ellipsoid for a perfect polytope $D^8_{11}\in Del(\Z^8,\cQ^8_{11})$, whose vertex set ($|\vertex D^8_{11}|=\textrm{44}$) is given below.
\begin{center}
\setlength{\extrarowheight}{4pt}
\begin{tabular}{|c|c|c|}
\hline$x_{8}=-1$ &$x_{8}=0$ &$x_{8}=1$\\
\hline\hline$\mathbf{[0^{6};1;-1]} \times 1$ &$\mathbf{[-1,0^{5};1;0]}\times 6$ &$\mathbf{[0^{7};1]} \times 1$\\
\cline{1-1}\cline{2-2}\cline{3-3}$\mathbf{[0,1^{4};0;-2;-1]}\times 5$ &$\mathbf{[0^{8}]} \times 1$ &$\mathbf{[0^{5};1;-1;1]} \times 1$\\
\cline{1-1}\cline{2-2}\cline{3-3}$\mathbf{[1^{6};-3;-1]} \times 1$ &$\mathbf{[0^{5},1;0^{2}]}\times 6$ &\\
\cline{1-1}\cline{2-2} &$\mathbf{[0^{4},1^{2};-1;0]}\times 15$ &\\
\cline{2-2} &$\mathbf{[0,1^{5};-3;0]}\times 6$ &\\
\cline{2-2} &$\mathbf{[1^{6};-4;0]} \times 1$ &\\
\hline
\end{tabular}
\end{center}
\begin{itemize}
\item $\Spec(D^8_{11})=\textrm{8, 9, 11, 12, 13, 15, 16, 17, 19, 20, 21, 23, 24, 27, 28, 31}$
\vspace{-2.5mm}\item $|Iso(D^8_{11},\cQ^8_{11})|=\textrm{240}$\vspace{-2.5mm}\item $l(D^8_{11})=\textrm{3}$
\vspace{-2.5mm}\item Antisymmetric
\vspace{-2.5mm}\item Maximally contained subpolytopes: $35-tope$, $H(2)$, $\frac12 H(5)$, $J(7,5)$
\end{itemize}
 
 Perfect affine quadratic lattice $\Aff(\Z^8,\cQ^8_{
12}[x-\cc])$, where $\cQ^8_{12}$ is given by
\begin{equation*}
\setlength{\extrarowheight}{1pt} \cQ^8_{12}(\x)=\x^t\; \begin{array}{|c|c|c|c|c|c|c|c|}
\hline  12 & 9 & 9 & 9 & 9 & -15 & 7 & 6\\
\hline  9 & 12 & 9 & 9 & 9 & -15 & 10 & 6\\
\hline  9 & 9 & 12 & 9 & 9 & -15 & 8 & 5\\
\hline  9 & 9 & 9 & 12 & 9 & -15 & 8 & 5\\
\hline  9 & 9 & 9 & 9 & 12 & -15 & 8 & 5\\
\hline  -15 & -15 & -15 & -15 & -15 & 24 & -12 & -9\\
\hline  7 & 10 & 8 & 8 & 8 & -12 & 12 & 5\\
\hline  6 & 6 & 5 & 5 & 5 & -9 & 5 & 6\\
\hline\end{array}
\;\x \quad \hbox{and}  \quad
\cc=\frac{1}{31}\;\begin{array}{|r|}
\hline 9\\
\hline 3\\
\hline 6\\
\hline 6\\
\hline 6\\
\hline 9\\
\hline 6\\
\hline -3\\
\hline\end{array}\;.
\end{equation*}
\begin{itemize}
\item $|O(\Z^8,\cQ^8_{12})|=\textrm{96}$; $S_{3}<O(\Z^8,\cQ^8_{12})$
\vspace{-2.5mm} \item $s(\Z^8,\cQ^8_{12})=\textrm{14}$
\vspace{-2.5mm} \item $\dim \QuadInv[\Z^8,\cQ^8_{12}]=\textrm{8}$
\end{itemize}
Inequality $\cQ_{12}^8[\x-\cc] \le \frac{126}{31}$ defines the Delaunay ellipsoid for a perfect polytope $D^8_{12}\in Del(\Z^8,\cQ^8_{12})$, whose vertex set ($|\vertex D^8_{12}|=\textrm{45}$) is given below.
\begin{center}
\setlength{\extrarowheight}{4pt}
\begin{tabular}{|c|c|c|}
\hline$x_{8}=-1$ &$x_{8}=0$ &$x_{8}=1$\\
\hline\hline$\mathbf{[0;-1^{4};-3;1;-1]} \times 1$ &$\mathbf{[-1,0^{4};-1;0^{2}]}\times 5$ &$\mathbf{[0^{7};1]} \times 1$\\
\cline{1-1}\cline{2-2}\cline{3-3}$\mathbf{[0;-1;-1^{2},0;-2;1;-1]}\times 3$ &$\mathbf{[0^{8}]} \times 1$ &$\mathbf{[0^{2};1^{3};2;0;1]} \times 1$\\
\cline{1-1}\cline{2-2}\cline{3-3}$\mathbf{[0;-1;0^{3};-1;1;-1]} \times 1$ &$\mathbf{[0^{6};1;0]} \times 1$ &$\mathbf{[0;1^{4};2;-1;1]} \times 1$\\
\cline{1-1}\cline{2-2}\cline{3-3}$\mathbf{[0^{5};-1;0;-1]} \times 1$ &$\mathbf{[0^{4},1;0^{3}]}\times 5$ &$\mathbf{[1^{5};3;0;1]} \times 1$\\
\cline{1-1}\cline{2-2}\cline{3-3}$\mathbf{[0;1;0^{5};-1]} \times 1$ &$\mathbf{[0^{3},1^{2};1;0^{2}]}\times 10$ &\\
\cline{1-1}\cline{2-2}$\mathbf{[1;0^{6};-1]} \times 1$ &$\mathbf{[0,1^{4};2;0^{2}]}\times 5$ &\\
\cline{1-1}\cline{2-2} &$\mathbf{[1^{5};3;0^{2}]} \times 1$ &\\
\cline{2-2} &$\mathbf{[0;-1^{2};0^{2};-1;1;0]} \times 1$ &\\
\cline{2-2} &$\mathbf{[0;-1;0;-1;0;-1;1;0]} \times 1$ &\\
\cline{2-2} &$\mathbf{[0;-1;0^{2};-1^{2};1;0]} \times 1$ &\\
\cline{2-2} &$\mathbf{[1;-1;0^{4};1;0]} \times 1$ &\\
\cline{2-2} &$\mathbf{[1;2;1^{3};3;-1;0]} \times 1$ &\\
\hline
\end{tabular}
\end{center}
\begin{itemize}
\item $\Spec(D^8_{12})=\textrm{2, 3, 4, 5, 6, 7}$
\vspace{-2.5mm}\item $|Iso(D^8_{12},\cQ^8_{12})|=\textrm{48}$\vspace{-2.5mm}\item $l(D^8_{12})=\textrm{3}$
\vspace{-2.5mm}\item Antisymmetric
\vspace{-2.5mm}\item Maximally contained subpolytopes: $G_6$, $H(2)$, $\frac12 H(5)$, $J(6,4)$
\end{itemize}

 Perfect affine quadratic lattice $\Aff(\Z^8,\cQ^8_{
13}[x-\cc])$, where $\cQ^8_{13}$ is given by
\begin{equation*}
\setlength{\extrarowheight}{1pt} \cQ^8_{13}(\x)=\x^t\; \begin{array}{|c|c|c|c|c|c|c|c|}
\hline  24 & 18 & 18 & 18 & 18 & 18 & 27 & 6\\
\hline  18 & 24 & 18 & 18 & 18 & 18 & 27 & 10\\
\hline  18 & 18 & 24 & 18 & 18 & 18 & 27 & 10\\
\hline  18 & 18 & 18 & 24 & 18 & 18 & 27 & 10\\
\hline  18 & 18 & 18 & 18 & 24 & 18 & 27 & 10\\
\hline  18 & 18 & 18 & 18 & 18 & 24 & 27 & 12\\
\hline  27 & 27 & 27 & 27 & 27 & 27 & 39 & 12\\
\hline  6 & 10 & 10 & 10 & 10 & 12 & 12 & 16\\
\hline\end{array}
\;\x \quad \hbox{and}  \quad
\cc=\frac{1}{74}\;\begin{array}{|r|}
\hline 22\\
\hline 16\\
\hline 16\\
\hline 16\\
\hline 16\\
\hline 13\\
\hline -40\\
\hline 9\\
\hline\end{array}\;.
\end{equation*}
\begin{itemize}
\item $|O(\Z^8,\cQ^8_{13})|=\textrm{96}$; $S_{4}<O(\Z^8,\cQ^8_{13})$
\vspace{-2.5mm} \item $s(\Z^8,\cQ^8_{13})=\textrm{2}$
\vspace{-2.5mm} \item $\dim \QuadInv[\Z^8,\cQ^8_{13}]=\textrm{12}$
\end{itemize}
Inequality $\cQ_{13}^8[\x-\cc] \le \frac{300}{37}$ defines the Delaunay ellipsoid for a perfect polytope $D^8_{13}\in Del(\Z^8,\cQ^8_{13})$, whose vertex set ($|\vertex D^8_{13}|=\textrm{44}$) is given below.
\begin{center}
\setlength{\extrarowheight}{4pt}
\begin{tabular}{|c|c|c|}
\hline$x_{8}=-1$ &$x_{8}=0$ &$x_{8}=1$\\
\hline\hline$\mathbf{[0;1^{5};-3;-1]} \times 1$ &$\mathbf{[-1,0^{5};1;0]}\times 6$ &$\mathbf{[0;-1^{3},0;-1;3;1]}\times 4$\\
\cline{1-1}\cline{2-2}\cline{3-3}$\mathbf{[1^{5};2;-4;-1]} \times 1$ &$\mathbf{[0^{8}]} \times 1$ &$\mathbf{[0^{5};-1;1^{2}]} \times 1$\\
\cline{1-1}\cline{2-2}\cline{3-3} &$\mathbf{[0^{5},1;0^{2}]}\times 6$ &$\mathbf{[0^{7};1]} \times 1$\\
\cline{2-2}\cline{3-3} &$\mathbf{[0^{4},1^{2};-1;0]}\times 15$ &$\mathbf{[1;0^{4};-1;0;1]} \times 1$\\
\cline{2-2}\cline{3-3} &$\mathbf{[0,1^{5};-3;0]}\times 6$ &\\
\cline{2-2} &$\mathbf{[1^{6};-4;0]} \times 1$ &\\
\hline
\end{tabular}
\end{center}
\begin{itemize}
\item $\Spec(D^8_{13})=\textrm{5, 8, 9, 11, 12, 13, 15, 16, 17, 19, 20, 21, 23, 24, 25, 27, 28}$
\vspace{-2.5mm}\item $|Iso(D^8_{13},\cQ^8_{13})|=\textrm{48}$\vspace{-2.5mm}\item $l(D^8_{13})=\textrm{3}$
\vspace{-2.5mm}\item Antisymmetric
\vspace{-2.5mm}\item Maximally contained subpolytopes: $35-tope$, $H(2)$, $\frac12 H(5)$, $J(7,5)$
\end{itemize}

 Perfect affine quadratic lattice $\Aff(\Z^8,\cQ^8_{
14}[x-\cc])$, where $\cQ^8_{14}$ is given by
\begin{equation*}
\setlength{\extrarowheight}{1pt} \cQ^8_{14}(\x)=\x^t\; \begin{array}{|c|c|c|c|c|c|c|c|}
\hline  9 & 5 & 5 & 5 & -15 & 0 & 3 & -12\\
\hline  5 & 9 & 5 & 5 & -15 & 0 & 3 & -12\\
\hline  5 & 5 & 9 & 5 & -15 & 0 & 3 & -12\\
\hline  5 & 5 & 5 & 9 & -15 & 0 & 3 & -12\\
\hline  -15 & -15 & -15 & -15 & 42 & 3 & -9 & 33\\
\hline  0 & 0 & 0 & 0 & 3 & 6 & 0 & 3\\
\hline  3 & 3 & 3 & 3 & -9 & 0 & 4 & -8\\
\hline  -12 & -12 & -12 & -12 & 33 & 3 & -8 & 28\\
\hline\end{array}
\;\x \quad \hbox{and}  \quad
\cc=\frac{1}{16}\;\begin{array}{|r|}
\hline 7\\
\hline 7\\
\hline 7\\
\hline 7\\
\hline 4\\
\hline 4\\
\hline 4\\
\hline 4\\
\hline\end{array}\;.
\end{equation*}
\begin{itemize}
\item $|O(\Z^8,\cQ^8_{14})|=\textrm{576}$; $S_{4}<O(\Z^8,\cQ^8_{14})$
\vspace{-2.5mm} \item $s(\Z^8,\cQ^8_{14})=\textrm{6}$
\vspace{-2.5mm} \item $\dim \QuadInv[\Z^8,\cQ^8_{14}]=\textrm{6}$
\end{itemize}
Inequality $\cQ_{14}^8[\x-\cc] \le \frac{41}{8}$ defines the Delaunay ellipsoid for a perfect polytope $D^8_{14}\in Del(\Z^8,\cQ^8_{14})$, whose vertex set ($|\vertex D^8_{14}|=\textrm{46}$) is given below.
\begin{center}
\setlength{\extrarowheight}{4pt}
\begin{tabular}{|c|c|c|}
\hline$x_{6}=-1$ &$x_{6}=0$ &$x_{6}=1$\\
\hline\hline$\mathbf{[1^{4};0;-1;1;2]} \times 1$ &$\mathbf{[0^{4};-1;0^{2};1]} \times 1$ &$\mathbf{[-1,0^{3};0;1;-1^{2}]}\times 4$\\
\cline{1-1}\cline{2-2}\cline{3-3}$\mathbf{[1^{4};0;-1;2^{2}]} \times 1$ &$\mathbf{[0^{8}]} \times 1$ &$\mathbf{[0^{4};-1;1;0;1]} \times 1$\\
\cline{1-1}\cline{2-2}\cline{3-3}$\mathbf{[1^{5};-1;2;1]} \times 1$ &$\mathbf{[0^{6};1;0]} \times 1$ &$\mathbf{[0^{5};1;-2;-1]} \times 1$\\
\cline{1-1}\cline{2-2}\cline{3-3} &$\mathbf{[0^{3},1;0^{4}]}\times 4$ &$\mathbf{[0^{5};1;-1^{2}]} \times 1$\\
\cline{2-2}\cline{3-3} &$\mathbf{[0^{2},1^{2};0^{2};1^{2}]}\times 6$ &$\mathbf{[0^{5};1;0^{2}]} \times 1$\\
\cline{2-2}\cline{3-3} &$\mathbf{[0,1^{3};0^{3};1]}\times 4$ &$\mathbf{[0^{5};1^{2};0]} \times 1$\\
\cline{2-2}\cline{3-3} &$\mathbf{[0,1^{3};1;0^{3}]}\times 4$ &$\mathbf{[0^{4};1^{2};-1;-2]} \times 1$\\
\cline{2-2}\cline{3-3} &$\mathbf{[0,1^{3};1;0;1;0]}\times 4$ &$\mathbf{[0^{3},1;0;1;0^{2}]}\times 4$\\
\cline{2-2}\cline{3-3} &$\mathbf{[1^{4};0^{2};1;2]} \times 1$ &\\
\cline{2-2} &$\mathbf{[1^{4};0^{2};2^{2}]} \times 1$ &\\
\cline{2-2} &$\mathbf{[1^{5};0^{3}]} \times 1$ &\\
\cline{2-2} &$\mathbf{[1^{5};0;2;1]} \times 1$ &\\
\hline
\end{tabular}
\end{center}
\begin{itemize}
\item $\Spec(D^8_{14})=\textrm{4, 6, 7, 8, 9, 10, 12, 13, 14, 15, 16, 18}$
\vspace{-2.5mm}\item $|Iso(D^8_{14},\cQ^8_{14})|=\textrm{288}$\vspace{-2.5mm}\item $l(D^8_{14})=\textrm{3}$
\vspace{-2.5mm}\item Antisymmetric
\vspace{-2.5mm}\item Maximally contained subpolytopes: $H(3)$, $\frac12 H(5)$, $J(7,5)$
\end{itemize}

 Perfect affine quadratic lattice $\Aff(\Z^8,\cQ^8_{
15}[x-\cc])$, where $\cQ^8_{15}$ is given by
\begin{equation*}
\setlength{\extrarowheight}{1pt} \cQ^8_{15}(\x)=\x^t\; \begin{array}{|c|c|c|c|c|c|c|c|}
\hline  8 & 6 & 6 & 6 & 6 & 6 & 9 & 2\\
\hline  6 & 8 & 6 & 6 & 6 & 6 & 9 & 3\\
\hline  6 & 6 & 8 & 6 & 6 & 6 & 9 & 3\\
\hline  6 & 6 & 6 & 8 & 6 & 6 & 9 & 3\\
\hline  6 & 6 & 6 & 6 & 8 & 6 & 9 & 3\\
\hline  6 & 6 & 6 & 6 & 6 & 8 & 9 & 4\\
\hline  9 & 9 & 9 & 9 & 9 & 9 & 13 & 4\\
\hline  2 & 3 & 3 & 3 & 3 & 4 & 4 & 4\\
\hline\end{array}
\;\x \quad \hbox{and}  \quad
\cc=\frac{1}{20}\;\begin{array}{|r|}
\hline 6\\
\hline 7\\
\hline 7\\
\hline 7\\
\hline 7\\
\hline 8\\
\hline -20\\
\hline -2\\
\hline\end{array}\;.
\end{equation*}
\begin{itemize}
\item $|O(\Z^8,\cQ^8_{15})|=\textrm{384}$; $S_{4}<O(\Z^8,\cQ^8_{15})$
\vspace{-2.5mm} \item $s(\Z^8,\cQ^8_{15})=\textrm{24}$
\vspace{-2.5mm} \item $\dim \QuadInv[\Z^8,\cQ^8_{15}]=\textrm{6}$
\end{itemize}
Inequality $\cQ_{15}^8[\x-\cc] \le \frac{27}{10}$ defines the Delaunay ellipsoid for a perfect polytope $D^8_{15}\in Del(\Z^8,\cQ^8_{15})$, whose vertex set ($|\vertex D^8_{15}|=\textrm{45}$) is given below.
\begin{center}
\setlength{\extrarowheight}{4pt}
\begin{tabular}{|c|c|c|}
\hline$x_{8}=-1$ &$x_{8}=0$ &$x_{8}=1$\\
\hline\hline$\mathbf{[0^{5};1;0;-1]} \times 1$ &$\mathbf{[-1,0^{5};1;0]}\times 6$ &$\mathbf{[0^{5};-1;1^{2}]} \times 1$\\
\cline{1-1}\cline{2-2}\cline{3-3}$\mathbf{[0;0,1^{3};1;-2;-1]}\times 4$ &$\mathbf{[0^{8}]} \times 1$ &$\mathbf{[0^{7};1]} \times 1$\\
\cline{1-1}\cline{2-2}\cline{3-3}$\mathbf{[0;1^{5};-3;-1]} \times 1$ &$\mathbf{[0^{5},1;0^{2}]}\times 6$ &$\mathbf{[1;0^{4};-1;0;1]} \times 1$\\
\cline{1-1}\cline{2-2}\cline{3-3}$\mathbf{[1^{5};2;-4;-1]} \times 1$ &$\mathbf{[0^{4},1^{2};-1;0]}\times 15$ &\\
\cline{1-1}\cline{2-2} &$\mathbf{[0,1^{5};-3;0]}\times 6$ &\\
\cline{2-2} &$\mathbf{[1^{6};-4;0]} \times 1$ &\\
\hline
\end{tabular}
\end{center}
\begin{itemize}
\item $\Spec(D^8_{15})=\textrm{3, 4, 5, 6, 7, 8, 9}$
\vspace{-2.5mm}\item $|Iso(D^8_{15},\cQ^8_{15})|=\textrm{192}$\vspace{-2.5mm}\item $l(D^8_{15})=\textrm{3}$
\vspace{-2.5mm}\item Antisymmetric
\vspace{-2.5mm}\item Maximally contained subpolytopes: $35-tope$, $H(2)$, $\frac12 H(5)$, $J(7,5)$
\end{itemize}

 Perfect affine quadratic lattice $\Aff(\Z^8,\cQ^8_{
16}[x-\cc])$, where $\cQ^8_{16}$ is given by
\begin{equation*}
\setlength{\extrarowheight}{1pt} \cQ^8_{16}(\x)=\x^t\; \begin{array}{|c|c|c|c|c|c|c|c|}
\hline  8 & 6 & 6 & 6 & 6 & 6 & 9 & 6\\
\hline  6 & 8 & 6 & 6 & 6 & 6 & 9 & 5\\
\hline  6 & 6 & 8 & 6 & 6 & 6 & 9 & 5\\
\hline  6 & 6 & 6 & 8 & 6 & 6 & 9 & 6\\
\hline  6 & 6 & 6 & 6 & 8 & 6 & 9 & 6\\
\hline  6 & 6 & 6 & 6 & 6 & 8 & 9 & 5\\
\hline  9 & 9 & 9 & 9 & 9 & 9 & 13 & 8\\
\hline  6 & 5 & 5 & 6 & 6 & 5 & 8 & 7\\
\hline\end{array}
\;\x \quad \hbox{and}  \quad
\cc=\frac{1}{14}\;\begin{array}{|r|}
\hline 4\\
\hline 5\\
\hline 5\\
\hline 4\\
\hline 4\\
\hline 5\\
\hline -14\\
\hline 2\\
\hline\end{array}\;.
\end{equation*}
\begin{itemize}
\item $|O(\Z^8,\cQ^8_{16})|=\textrm{288}$; $S_{3}<O(\Z^8,\cQ^8_{16})$
\vspace{-2.5mm} \item $s(\Z^8,\cQ^8_{16})=\textrm{24}$
\vspace{-2.5mm} \item $\dim \QuadInv[\Z^8,\cQ^8_{16}]=\textrm{5}$
\end{itemize}
Inequality $\cQ_{16}^8[\x-\cc] \le \frac{19}{7}$ defines the Delaunay ellipsoid for a perfect polytope $D^8_{16}\in Del(\Z^8,\cQ^8_{16})$, whose vertex set ($|\vertex D^8_{16}|=\textrm{45}$) is given below.
\begin{center}
\setlength{\extrarowheight}{4pt}
\begin{tabular}{|c|c|c|}
\hline$x_{8}=-1$ &$x_{8}=0$ &$x_{8}=1$\\
\hline\hline$\mathbf{[0^{6};1;-1]} \times 1$ &$\mathbf{[-1,0^{5};1;0]}\times 6$ &$\mathbf{[0^{7};1]} \times 1$\\
\cline{1-1}\cline{2-2}\cline{3-3}$\mathbf{[1;0^{2};1^{2};0;-1^{2}]} \times 1$ &$\mathbf{[0^{8}]} \times 1$ &$\mathbf{[0^{5};1;-1;1]} \times 1$\\
\cline{1-1}\cline{2-2}\cline{3-3} &$\mathbf{[0^{5},1;0^{2}]}\times 6$ &$\mathbf{[0;0,1;0^{3};-1;1]}\times 2$\\
\cline{2-2}\cline{3-3} &$\mathbf{[0^{4},1^{2};-1;0]}\times 15$ &$\mathbf{[0;1^{2};0^{2};1;-2;1]} \times 1$\\
\cline{2-2}\cline{3-3} &$\mathbf{[0,1^{5};-3;0]}\times 6$ &$\mathbf{[0;1^{2};0;1^{2};-3;1]} \times 1$\\
\cline{2-2}\cline{3-3} &$\mathbf{[1^{6};-4;0]} \times 1$ &$\mathbf{[0;1^{3};0;1;-3;1]} \times 1$\\
\cline{2-2}\cline{3-3} & &$\mathbf{[1^{3};0^{2};1;-3;1]} \times 1$\\
\hline
\end{tabular}
\end{center}
\begin{itemize}
\item $\Spec(D^8_{16})=\textrm{3, 4, 5, 6, 7, 8, 9}$
\vspace{-2.5mm}\item $|Iso(D^8_{16},\cQ^8_{16})|=\textrm{144}$\vspace{-2.5mm}\item $l(D^8_{16})=\textrm{3}$
\vspace{-2.5mm}\item Antisymmetric
\vspace{-2.5mm}\item Maximally contained subpolytopes: $35-tope$, $H(2)$, $\frac12 H(5)$, $J(7,5)$
\end{itemize}

 Perfect affine quadratic lattice $\Aff(\Z^8,\cQ^8_{
17}[x-\cc])$, where $\cQ^8_{17}$ is given by
\begin{equation*}
\setlength{\extrarowheight}{1pt} \cQ^8_{17}(\x)=\x^t\; \begin{array}{|c|c|c|c|c|c|c|c|}
\hline  8 & 6 & 6 & 6 & 6 & 6 & 9 & 4\\
\hline  6 & 8 & 6 & 6 & 6 & 6 & 9 & 4\\
\hline  6 & 6 & 8 & 6 & 6 & 6 & 9 & 4\\
\hline  6 & 6 & 6 & 8 & 6 & 6 & 9 & 6\\
\hline  6 & 6 & 6 & 6 & 8 & 6 & 9 & 4\\
\hline  6 & 6 & 6 & 6 & 6 & 8 & 9 & 4\\
\hline  9 & 9 & 9 & 9 & 9 & 9 & 13 & 6\\
\hline  4 & 4 & 4 & 6 & 4 & 4 & 6 & 6\\
\hline\end{array}
\;\x \quad \hbox{and}  \quad
\cc=\frac{1}{10}\;\begin{array}{|r|}
\hline 3\\
\hline 3\\
\hline 3\\
\hline 2\\
\hline 3\\
\hline 3\\
\hline -8\\
\hline 1\\
\hline\end{array}\;.
\end{equation*}
\begin{itemize}
\item $|O(\Z^8,\cQ^8_{17})|=\textrm{5760}$; $S_{6}<O(\Z^8,\cQ^8_{17})$
\vspace{-2.5mm} \item $s(\Z^8,\cQ^8_{17})=\textrm{2}$
\vspace{-2.5mm} \item $\dim \QuadInv[\Z^8,\cQ^8_{17}]=\textrm{4}$
\end{itemize}
Inequality $\cQ_{17}^8[\x-\cc] \le \frac{27}{10}$ defines the Delaunay ellipsoid for a perfect polytope $D^8_{17}\in Del(\Z^8,\cQ^8_{17})$, whose vertex set ($|\vertex D^8_{17}|=\textrm{44}$) is given below.
\begin{center}
\setlength{\extrarowheight}{4pt}
\begin{tabular}{|c|c|c|}
\hline$x_{8}=-1$ &$x_{8}=0$ &$x_{8}=1$\\
\hline\hline$\mathbf{[0^{3};1;0^{3};-1]} \times 1$ &$\mathbf{[-1,0^{5};1;0]}\times 6$ &$\mathbf{[0^{3};-1;0^{2};1^{2}]} \times 1$\\
\cline{1-1}\cline{2-2}\cline{3-3}$\mathbf{[1^{3};2;1^{2};-4;-1]} \times 1$ &$\mathbf{[0^{8}]} \times 1$ &$\mathbf{[0^{3};-1;0;1;0;1]} \times 1$\\
\cline{1-1}\cline{2-2}\cline{3-3} &$\mathbf{[0^{5},1;0^{2}]}\times 6$ &$\mathbf{[0^{3};-1;1;0^{2};1]} \times 1$\\
\cline{2-2}\cline{3-3} &$\mathbf{[0^{4},1^{2};-1;0]}\times 15$ &$\mathbf{[0^{7};1]} \times 1$\\
\cline{2-2}\cline{3-3} &$\mathbf{[0,1^{5};-3;0]}\times 6$ &$\mathbf{[0^{2},1;-1;0^{3};1]}\times 3$\\
\cline{2-2}\cline{3-3} &$\mathbf{[1^{6};-4;0]} \times 1$ &\\
\hline
\end{tabular}
\end{center}
\begin{itemize}
\item $\Spec(D^8_{17})=\textrm{2, 3, 4, 5, 6, 7, 8, 9}$
\vspace{-2.5mm}\item $|Iso(D^8_{17},\cQ^8_{17})|=\textrm{2880}$\vspace{-2.5mm}\item $l(D^8_{17})=\textrm{3}$
\vspace{-2.5mm}\item Antisymmetric
\vspace{-2.5mm}\item Maximally contained subpolytopes: $35-tope$, $H(2)$, $\frac12 H(5)$, $J(8,6)$
\end{itemize}

 Perfect affine quadratic lattice $\Aff(\Z^8,\cQ^8_{
18}[x-\cc])$, where $\cQ^8_{18}$ is given by
\begin{equation*}
\setlength{\extrarowheight}{1pt} \cQ^8_{18}(\x)=\x^t\; \begin{array}{|c|c|c|c|c|c|c|c|}
\hline  11 & 6 & 6 & 6 & -7 & 11 & 4 & 3\\
\hline  6 & 11 & 6 & 6 & -7 & 11 & 4 & 3\\
\hline  6 & 6 & 11 & 6 & -7 & 11 & 4 & 3\\
\hline  6 & 6 & 6 & 11 & -7 & 11 & 4 & 3\\
\hline  -7 & -7 & -7 & -7 & 12 & -8 & -4 & 4\\
\hline  11 & 11 & 11 & 11 & -8 & 20 & 4 & 8\\
\hline  4 & 4 & 4 & 4 & -4 & 4 & 8 & 0\\
\hline  3 & 3 & 3 & 3 & 4 & 8 & 0 & 13\\
\hline\end{array}
\;\x \quad \hbox{and}  \quad
\cc=\frac{1}{76}\;\begin{array}{|r|}
\hline 6\\
\hline 6\\
\hline 6\\
\hline 6\\
\hline 23\\
\hline 25\\
\hline 25\\
\hline 10\\
\hline\end{array}\;.
\end{equation*}
\begin{itemize}
\item $|O(\Z^8,\cQ^8_{18})|=\textrm{288}$; $S_{4}<O(\Z^8,\cQ^8_{18})$
\vspace{-2.5mm} \item $s(\Z^8,\cQ^8_{18})=\textrm{6}$
\vspace{-2.5mm} \item $\dim \QuadInv[\Z^8,\cQ^8_{18}]=\textrm{8}$
\end{itemize}
Inequality $\cQ_{18}^8[\x-\cc] \le \frac{501}{76}$ defines the Delaunay ellipsoid for a perfect polytope $D^8_{18}\in Del(\Z^8,\cQ^8_{18})$, whose vertex set ($|\vertex D^8_{18}|=\textrm{44}$) is given below.
\begin{center}
\setlength{\extrarowheight}{4pt}
\begin{tabular}{|c|c|c|c|}
\hline$x_{6}=-1$ &$x_{6}=0$ &$x_{6}=1$ &$x_{6}=2$\\
\hline\hline$\mathbf{[0,1^{3};1;-1;0^{2}]}\times 4$ &$\mathbf{[0^{4};-1;0^{2};1]} \times 1$ &$\mathbf{[-1^{4};-2;1^{3}]} \times 1$ &$\mathbf{[-1^{5};2;1;0]} \times 1$\\
\cline{1-1}\cline{2-2}\cline{3-3}\cline{4-4}$\mathbf{[1^{5};-1^{2};0]} \times 1$ &$\mathbf{[0^{8}]} \times 1$ &$\mathbf{[-1^{3},0;-1;1^{3}]}\times 4$ &$\mathbf{[-1^{5};2;1^{2}]} \times 1$\\
\cline{1-1}\cline{2-2}\cline{3-3}\cline{4-4}$\mathbf{[1^{5};-1;0^{2}]} \times 1$ &$\mathbf{[0^{7};1]} \times 1$ &$\mathbf{[-1,0^{3};0;1;0^{2}]}\times 4$ &\\
\cline{1-1}\cline{2-2}\cline{3-3}$\mathbf{[1^{4};2;-1;0^{2}]} \times 1$ &$\mathbf{[0^{6};1;0]} \times 1$ &$\mathbf{[-1,0^{3};0;1^{2};0]}\times 4$ &\\
\cline{1-1}\cline{2-2}\cline{3-3} &$\mathbf{[0^{6};1^{2}]} \times 1$ &$\mathbf{[0^{5};1;0^{2}]} \times 1$ &\\
\cline{2-2}\cline{3-3} &$\mathbf{[0^{3},1;0^{4}]}\times 4$ &$\mathbf{[0^{4};1^{2};0^{2}]} \times 1$ &\\
\cline{2-2}\cline{3-3} &$\mathbf{[0^{2},1^{2};1;0^{3}]}\times 6$ &$\mathbf{[0^{4};1^{3};0]} \times 1$ &\\
\cline{2-2}\cline{3-3} &$\mathbf{[0,1^{3};2;0^{2};-1]}\times 4$ & &\\
\hline
\end{tabular}
\end{center}
\begin{itemize}
\item $\Spec(D^8_{18})=\textrm{5, 8, 9, 10, 11, 12, 13, 14, 15, 16, 17, 18, 19, 20, 21, 24}$
\vspace{-2.5mm}\item $|Iso(D^8_{18},\cQ^8_{18})|=\textrm{144}$\vspace{-2.5mm}\item $l(D^8_{18})=\textrm{3}$
\vspace{-2.5mm}\item Antisymmetric
\vspace{-2.5mm}\item Maximally contained subpolytopes: $H(2)$, $\frac12 H(5)$, $J(7,5)$
\end{itemize}

 Perfect affine quadratic lattice $\Aff(\Z^8,\cQ^8_{
19}[x-\cc])$, where $\cQ^8_{19}$ is given by
\begin{equation*}
\setlength{\extrarowheight}{1pt} \cQ^8_{19}(\x)=\x^t\; \begin{array}{|c|c|c|c|c|c|c|c|}
\hline  8 & 5 & 5 & 5 & 1 & 1 & 8 & 4\\
\hline  5 & 8 & 5 & 5 & 1 & 1 & 8 & 4\\
\hline  5 & 5 & 8 & 5 & 1 & 1 & 8 & 4\\
\hline  5 & 5 & 5 & 8 & 1 & 1 & 8 & 4\\
\hline  1 & 1 & 1 & 1 & 7 & 3 & 2 & 0\\
\hline  1 & 1 & 1 & 1 & 3 & 7 & 6 & 4\\
\hline  8 & 8 & 8 & 8 & 2 & 6 & 15 & 8\\
\hline  4 & 4 & 4 & 4 & 0 & 4 & 8 & 8\\
\hline\end{array}
\;\x \quad \hbox{and}  \quad
\cc=\frac{1}{46}\;\begin{array}{|r|}
\hline -11\\
\hline -11\\
\hline -11\\
\hline -11\\
\hline 29\\
\hline -32\\
\hline 49\\
\hline 12\\
\hline\end{array}\;.
\end{equation*}
\begin{itemize}
\item $|O(\Z^8,\cQ^8_{19})|=\textrm{576}$; $S_{4}<O(\Z^8,\cQ^8_{19})$
\vspace{-2.5mm} \item $s(\Z^8,\cQ^8_{19})=\textrm{8}$
\vspace{-2.5mm} \item $\dim \QuadInv[\Z^8,\cQ^8_{19}]=\textrm{6}$
\end{itemize}
Inequality $\cQ_{19}^8[\x-\cc] \le \frac{229}{46}$ defines the Delaunay ellipsoid for a perfect polytope $D^8_{19}\in Del(\Z^8,\cQ^8_{19})$, whose vertex set ($|\vertex D^8_{19}|=\textrm{49}$) is given below.
\begin{center}
\setlength{\extrarowheight}{4pt}
\begin{tabular}{|c|c|c|}
\hline$x_{8}=-1$ &$x_{8}=0$ &$x_{8}=1$\\
\hline\hline$\mathbf{[0^{6};1;-1]} \times 1$ &$\mathbf{[-1^{4};1;-2;3;0]} \times 1$ &$\mathbf{[-1^{4};1;-3;3;1]} \times 1$\\
\cline{1-1}\cline{2-2}\cline{3-3} &$\mathbf{[-1^{3},0;1;-2;3;0]}\times 4$ &$\mathbf{[-1^{4};1;-2;3;1]} \times 1$\\
\cline{2-2}\cline{3-3} &$\mathbf{[-1^{2},0^{2};1;-1;2;0]}\times 6$ &$\mathbf{[-1^{4};2;-3;3;1]} \times 1$\\
\cline{2-2}\cline{3-3} &$\mathbf{[-1,0^{3};0^{2};1;0]}\times 4$ &$\mathbf{[-1^{3},0;1;-2;2;1]}\times 4$\\
\cline{2-2}\cline{3-3} &$\mathbf{[-1,0^{3};1;-1;1;0]}\times 4$ &$\mathbf{[-1,0^{3};1;-1;1^{2}]}\times 4$\\
\cline{2-2}\cline{3-3} &$\mathbf{[0^{8}]} \times 1$ &$\mathbf{[0^{7};1]} \times 1$\\
\cline{2-2}\cline{3-3} &$\mathbf{[0^{6};1;0]} \times 1$ &$\mathbf{[0^{4};1;-1;0;1]} \times 1$\\
\cline{2-2}\cline{3-3} &$\mathbf{[0^{5};1;0^{2}]} \times 1$ &$\mathbf{[0^{4};1;0^{2};1]} \times 1$\\
\cline{2-2}\cline{3-3} &$\mathbf{[0^{4};1;-1;1;0]} \times 1$ &\\
\cline{2-2} &$\mathbf{[0^{4};1;0^{3}]} \times 1$ &\\
\cline{2-2} &$\mathbf{[0^{3},1;0^{4}]}\times 4$ &\\
\cline{2-2} &$\mathbf{[0^{2},1^{2};0;1;-1;0]}\times 6$ &\\
\hline
\end{tabular}
\end{center}
\begin{itemize}
\item $\Spec(D^8_{19})=\textrm{4, 6, 7, 8, 9, 10, 12, 13, 14, 15, 16, 18}$
\vspace{-2.5mm}\item $|Iso(D^8_{19},\cQ^8_{19})|=\textrm{288}$\vspace{-2.5mm}\item $l(D^8_{19})=\textrm{3}$
\vspace{-2.5mm}\item Antisymmetric
\vspace{-2.5mm}\item Maximally contained subpolytopes: $H(3)$, $\frac12 H(5)$, $J(7,5)$
\end{itemize}

 Perfect affine quadratic lattice $\Aff(\Z^8,\cQ^8_{
20}[x-\cc])$, where $\cQ^8_{20}$ is given by
\begin{equation*}
\setlength{\extrarowheight}{1pt} \cQ^8_{20}(\x)=\x^t\; \begin{array}{|c|c|c|c|c|c|c|c|}
\hline  12 & 9 & 9 & 9 & 9 & -15 & 3 & 3\\
\hline  9 & 12 & 9 & 9 & 9 & -15 & 4 & 5\\
\hline  9 & 9 & 12 & 9 & 9 & -15 & 4 & 3\\
\hline  9 & 9 & 9 & 12 & 9 & -15 & 4 & 5\\
\hline  9 & 9 & 9 & 9 & 12 & -15 & 4 & 4\\
\hline  -15 & -15 & -15 & -15 & -15 & 24 & -6 & -6\\
\hline  3 & 4 & 4 & 4 & 4 & -6 & 4 & 3\\
\hline  3 & 5 & 3 & 5 & 4 & -6 & 3 & 5\\
\hline\end{array}
\;\x \quad \hbox{and}  \quad
\cc=\frac{1}{20}\;\begin{array}{|r|}
\hline 6\\
\hline 7\\
\hline 7\\
\hline 7\\
\hline 7\\
\hline 13\\
\hline -3\\
\hline 0\\
\hline\end{array}\;.
\end{equation*}
\begin{itemize}
\item $|O(\Z^8,\cQ^8_{20})|=\textrm{48}$; $S_{3}<O(\Z^8,\cQ^8_{20})$
\vspace{-2.5mm} \item $s(\Z^8,\cQ^8_{20})=\textrm{4}$
\vspace{-2.5mm} \item $\dim \QuadInv[\Z^8,\cQ^8_{20}]=\textrm{11}$
\end{itemize}
Inequality $\cQ_{20}^8[\x-\cc] \le \frac{81}{20}$ defines the Delaunay ellipsoid for a perfect polytope $D^8_{20}\in Del(\Z^8,\cQ^8_{20})$, whose vertex set ($|\vertex D^8_{20}|=\textrm{47}$) is given below.
\begin{center}
\setlength{\extrarowheight}{4pt}
\begin{tabular}{|c|c|c|}
\hline$x_{8}=-1$ &$x_{8}=0$ &$x_{8}=1$\\
\hline\hline$\mathbf{[0^{3};1;0^{3};-1]} \times 1$ &$\mathbf{[-1,0^{4};-1;0^{2}]}\times 5$ &$\mathbf{[0;-1;0^{3};-1^{2};1]} \times 1$\\
\cline{1-1}\cline{2-2}\cline{3-3}$\mathbf{[0;1;0^{5};-1]} \times 1$ &$\mathbf{[0^{8}]} \times 1$ &$\mathbf{[0^{3};-1;0;-1^{2};1]} \times 1$\\
\cline{1-1}\cline{2-2}\cline{3-3}$\mathbf{[0;1;0;1;0;1^{2};-1]} \times 1$ &$\mathbf{[0^{6};1;0]} \times 1$ &$\mathbf{[0^{7};1]} \times 1$\\
\cline{1-1}\cline{2-2}\cline{3-3}$\mathbf{[0;1;0;1^{3};0;-1]} \times 1$ &$\mathbf{[0^{4},1;0^{3}]}\times 5$ &$\mathbf{[0^{2};1;0^{3};-1;1]} \times 1$\\
\cline{1-1}\cline{2-2}\cline{3-3}$\mathbf{[0;1^{4};2;0;-1]} \times 1$ &$\mathbf{[0^{3},1^{2};1;0^{2}]}\times 10$ &$\mathbf{[0^{2};1;0;1^{2};-1;1]} \times 1$\\
\cline{1-1}\cline{2-2}\cline{3-3}$\mathbf{[1^{2};0;1^{2};2;0;-1]} \times 1$ &$\mathbf{[0,1^{4};2;0^{2}]}\times 5$ &$\mathbf{[1;0;1;0^{2};1;-1;1]} \times 1$\\
\cline{1-1}\cline{2-2}\cline{3-3}$\mathbf{[1^{2};0;1^{2};2;1;-1]} \times 1$ &$\mathbf{[1^{5};3;0^{2}]} \times 1$ &$\mathbf{[1;0;1;0^{2};1;0;1]} \times 1$\\
\cline{1-1}\cline{2-2}\cline{3-3} &$\mathbf{[0^{4};1;0;-1;0]} \times 1$ &\\
\cline{2-2} &$\mathbf{[0^{3};1;0^{2};-1;0]} \times 1$ &\\
\cline{2-2} &$\mathbf{[0^{2};1;0^{3};-1;0]} \times 1$ &\\
\cline{2-2} &$\mathbf{[0;1;0^{4};-1;0]} \times 1$ &\\
\cline{2-2} &$\mathbf{[0;1^{4};2;-1;0]} \times 1$ &\\
\hline
\end{tabular}
\end{center}
\begin{itemize}
\item $\Spec(D^8_{20})=\textrm{3, 4, 5, 6, 7, 8, 9, 10, 11, 12, 13, 14}$
\vspace{-2.5mm}\item $|Iso(D^8_{20},\cQ^8_{20})|=\textrm{24}$\vspace{-2.5mm}\item $l(D^8_{20})=\textrm{3}$
\vspace{-2.5mm}\item Antisymmetric
\vspace{-2.5mm}\item Maximally contained subpolytopes: $G_6$, $H(3)$, $\frac12 H(5)$, $J(6,4)$
\end{itemize}

 Perfect affine quadratic lattice $\Aff(\Z^8,\cQ^8_{
21}[x-\cc])$, where $\cQ^8_{21}$ is given by
\begin{equation*}
\setlength{\extrarowheight}{1pt} \cQ^8_{21}(\x)=\x^t\; \begin{array}{|c|c|c|c|c|c|c|c|}
\hline  12 & 9 & 9 & 9 & 9 & -15 & 6 & 5\\
\hline  9 & 12 & 9 & 9 & 9 & -15 & 3 & 4\\
\hline  9 & 9 & 12 & 9 & 9 & -15 & 4 & 3\\
\hline  9 & 9 & 9 & 12 & 9 & -15 & 4 & 3\\
\hline  9 & 9 & 9 & 9 & 12 & -15 & 4 & 3\\
\hline  -15 & -15 & -15 & -15 & -15 & 24 & -6 & -5\\
\hline  6 & 3 & 4 & 4 & 4 & -6 & 6 & 4\\
\hline  5 & 4 & 3 & 3 & 3 & -5 & 4 & 5\\
\hline\end{array}
\;\x \quad \hbox{and}  \quad
\cc=\frac{1}{50}\;\begin{array}{|r|}
\hline 19\\
\hline 20\\
\hline 17\\
\hline 17\\
\hline 17\\
\hline 37\\
\hline 3\\
\hline -6\\
\hline\end{array}\;.
\end{equation*}
\begin{itemize}
\item $|O(\Z^8,\cQ^8_{21})|=\textrm{96}$; $S_{3}<O(\Z^8,\cQ^8_{21})$
\vspace{-2.5mm} \item $s(\Z^8,\cQ^8_{21})=\textrm{4}$
\vspace{-2.5mm} \item $\dim \QuadInv[\Z^8,\cQ^8_{21}]=\textrm{8}$
\end{itemize}
Inequality $\cQ_{21}^8[\x-\cc] \le \frac{201}{50}$ defines the Delaunay ellipsoid for a perfect polytope $D^8_{21}\in Del(\Z^8,\cQ^8_{21})$, whose vertex set ($|\vertex D^8_{21}|=\textrm{47}$) is given below.
\begin{center}
\setlength{\extrarowheight}{4pt}
\begin{tabular}{|c|c|c|}
\hline$x_{8}=-1$ &$x_{8}=0$ &$x_{8}=1$\\
\hline\hline$\mathbf{[0;1;0^{4};1;-1]} \times 1$ &$\mathbf{[-1,0^{4};-1;0^{2}]}\times 5$ &$\mathbf{[-1;0^{4};-1;0;1]} \times 1$\\
\cline{1-1}\cline{2-2}\cline{3-3}$\mathbf{[0;1;0^{2},1;1^{2};-1]}\times 3$ &$\mathbf{[0^{8}]} \times 1$ &$\mathbf{[0;-1;0^{3};-1^{2};1]} \times 1$\\
\cline{1-1}\cline{2-2}\cline{3-3}$\mathbf{[1;0^{6};-1]} \times 1$ &$\mathbf{[0^{6};1;0]} \times 1$ &$\mathbf{[0^{7};1]} \times 1$\\
\cline{1-1}\cline{2-2}\cline{3-3}$\mathbf{[1^{2};0^{3};1^{2};-1]} \times 1$ &$\mathbf{[0^{4},1;0^{3}]}\times 5$ &$\mathbf{[1;0;1^{3};2;-1;1]} \times 1$\\
\cline{1-1}\cline{2-2}\cline{3-3}$\mathbf{[1^{2};0,1^{2};2;0;-1]}\times 3$ &$\mathbf{[0^{3},1^{2};1;0^{2}]}\times 10$ &\\
\cline{1-1}\cline{2-2}$\mathbf{[2;1^{4};3;0;-1]} \times 1$ &$\mathbf{[0,1^{4};2;0^{2}]}\times 5$ &\\
\cline{1-1}\cline{2-2} &$\mathbf{[1^{5};3;0^{2}]} \times 1$ &\\
\cline{2-2} &$\mathbf{[-1;0^{4};-1;1;0]} \times 1$ &\\
\cline{2-2} &$\mathbf{[-1;1;0^{4};1;0]} \times 1$ &\\
\cline{2-2} &$\mathbf{[1;0^{5};-1;0]} \times 1$ &\\
\cline{2-2} &$\mathbf{[1;0;1^{3};2;-1;0]} \times 1$ &\\
\cline{2-2} &$\mathbf{[2;1^{4};3;-1;0]} \times 1$ &\\
\hline
\end{tabular}
\end{center}
\begin{itemize}
\item $\Spec(D^8_{21})=\textrm{3, 4, 5, 6, 7, 8, 9, 10, 11, 12, 13}$
\vspace{-2.5mm}\item $|Iso(D^8_{21},\cQ^8_{21})|=\textrm{48}$\vspace{-2.5mm}\item $l(D^8_{21})=\textrm{3}$
\vspace{-2.5mm}\item Antisymmetric
\vspace{-2.5mm}\item Maximally contained subpolytopes: $G_6$, $H(3)$, $\frac12 H(5)$, $J(6,4)$
\end{itemize}
 
 Perfect affine quadratic lattice $\Aff(\Z^8,\cQ^8_{
22}[x-\cc])$, where $\cQ^8_{22}$ is given by
\begin{equation*}
\setlength{\extrarowheight}{1pt} \cQ^8_{22}(\x)=\x^t\; \begin{array}{|c|c|c|c|c|c|c|c|}
\hline  3 & 1 & 0 & 0 & 0 & 0 & 0 & 0\\
\hline  1 & 4 & 2 & 2 & 2 & 2 & 2 & 2\\
\hline  0 & 2 & 2 & 1 & 1 & 1 & 1 & 1\\
\hline  0 & 2 & 1 & 2 & 1 & 1 & 1 & 1\\
\hline  0 & 2 & 1 & 1 & 2 & 1 & 1 & 1\\
\hline  0 & 2 & 1 & 1 & 1 & 2 & 1 & 1\\
\hline  0 & 2 & 1 & 1 & 1 & 1 & 2 & 1\\
\hline  0 & 2 & 1 & 1 & 1 & 1 & 1 & 2\\
\hline\end{array}
\;\x \quad \hbox{and}  \quad
\cc=\frac{1}{10}\;\begin{array}{|r|}
\hline 8\\
\hline -9\\
\hline 4\\
\hline 4\\
\hline 4\\
\hline 4\\
\hline 4\\
\hline 4\\
\hline\end{array}\;.
\end{equation*}
\begin{itemize}
\item $|O(\Z^8,\cQ^8_{22})|=\textrm{645120}$; $S_{7}<O(\Z^8,\cQ^8_{22})$
\vspace{-2.5mm} \item $s(\Z^8,\cQ^8_{22})=\textrm{84}$
\vspace{-2.5mm} \item $\dim \QuadInv[\Z^8,\cQ^8_{22}]=\textrm{2}$
\end{itemize}
Inequality $\cQ_{22}^8[\x-\cc] \le \frac{9}{5}$ defines the Delaunay ellipsoid for a perfect polytope $D^8_{22}\in Del(\Z^8,\cQ^8_{22})$, whose vertex set ($|\vertex D^8_{22}|=\textrm{79}$) is given below.
\begin{center}
\setlength{\extrarowheight}{4pt}
\begin{tabular}{|c|c|c|}
\hline$x_{1}=0$ &$x_{1}=1$ &$x_{1}=2$\\
\hline\hline$\mathbf{[0^{8}]} \times 1$ &$\mathbf{[1;-3;1^{6}]} \times 1$ &$\mathbf{[2;-3;1^{6}]} \times 1$\\
\cline{1-1}\cline{2-2}\cline{3-3}$\mathbf{[0^{2};0^{5},1]}\times 6$ &$\mathbf{[1;-2;0^{2},1^{4}]}\times 15$ &\\
\cline{1-1}\cline{2-2}$\mathbf{[0;1;-1,0^{5}]}\times 6$ &$\mathbf{[1;-2;0,1^{5}]}\times 6$ &\\
\cline{1-1}\cline{2-2}$\mathbf{[0;1;0^{6}]} \times 1$ &$\mathbf{[1;-1;0^{4},1^{2}]}\times 15$ &\\
\cline{1-1}\cline{2-2} &$\mathbf{[1;-1;0^{3},1^{3}]}\times 20$ &\\
\cline{2-2} &$\mathbf{[1;0^{7}]} \times 1$ &\\
\cline{2-2} &$\mathbf{[1;0;0^{5},1]}\times 6$ &\\
\hline
\end{tabular}
\end{center}
\begin{itemize}
\item $\Spec(D^8_{22})=\textrm{2, 3, 4, 5, 6}$
\vspace{-2.5mm}\item $|Iso(D^8_{22},\cQ^8_{22})|=\textrm{322560}$\vspace{-2.5mm}\item $l(D^8_{22})=\textrm{3}$
\vspace{-2.5mm}\item Antisymmetric
\vspace{-2.5mm}\item Maximally contained subpolytopes: $H(3)$, $\frac12 H(7)$, $J(8,6)$
\end{itemize}

 Perfect affine quadratic lattice $\Aff(\Z^8,\cQ^8_{
23}[x-\cc])$, where $\cQ^8_{23}$ is given by
\begin{equation*}
\setlength{\extrarowheight}{1pt} \cQ^8_{23}(\x)=\x^t\; \begin{array}{|c|c|c|c|c|c|c|c|}
\hline  7 & 4 & 4 & 4 & 4 & 2 & 7 & 4\\
\hline  4 & 7 & 4 & 4 & 4 & 2 & 7 & 4\\
\hline  4 & 4 & 7 & 4 & 4 & 2 & 7 & 4\\
\hline  4 & 4 & 4 & 7 & 4 & 2 & 7 & 4\\
\hline  4 & 4 & 4 & 4 & 7 & 2 & 7 & 4\\
\hline  2 & 2 & 2 & 2 & 2 & 7 & 7 & 0\\
\hline  7 & 7 & 7 & 7 & 7 & 7 & 14 & 4\\
\hline  4 & 4 & 4 & 4 & 4 & 0 & 4 & 8\\
\hline\end{array}
\;\x \quad \hbox{and}  \quad
\cc=\frac{1}{20}\;\begin{array}{|r|}
\hline -7\\
\hline -7\\
\hline -7\\
\hline -7\\
\hline -7\\
\hline -7\\
\hline 27\\
\hline 14\\
\hline\end{array}\;.
\end{equation*}
\begin{itemize}
\item $|O(\Z^8,\cQ^8_{23})|=\textrm{1920}$; $S_{5}<O(\Z^8,\cQ^8_{23})$
\vspace{-2.5mm} \item $s(\Z^8,\cQ^8_{23})=\textrm{10}$
\vspace{-2.5mm} \item $\dim \QuadInv[\Z^8,\cQ^8_{23}]=\textrm{4}$
\end{itemize}
Inequality $\cQ_{23}^8[\x-\cc] \le \frac{49}{10}$ defines the Delaunay ellipsoid for a perfect polytope $D^8_{23}\in Del(\Z^8,\cQ^8_{23})$, whose vertex set ($|\vertex D^8_{23}|=\textrm{49}$) is given below.
\begin{center}
\setlength{\extrarowheight}{4pt}
\begin{tabular}{|c|c|c|}
\hline$x_{8}=0$ &$x_{8}=1$ &$x_{8}=2$\\
\hline\hline$\mathbf{[-1,0^{5};1;0]}\times 6$ &$\mathbf{[-1^{5};-2;4;1]} \times 1$ &$\mathbf{[-1^{6};3;2]} \times 1$\\
\cline{1-1}\cline{2-2}\cline{3-3}$\mathbf{[0^{8}]} \times 1$ &$\mathbf{[-1^{6};3;1]} \times 1$ &\\
\cline{1-1}\cline{2-2}$\mathbf{[0^{6};1;0]} \times 1$ &$\mathbf{[-1^{4},0;-1;3;1]}\times 5$ &\\
\cline{1-1}\cline{2-2}$\mathbf{[0^{5},1;0^{2}]}\times 6$ &$\mathbf{[-1^{3},0^{2};-1;2;1]}\times 10$ &\\
\cline{1-1}\cline{2-2} &$\mathbf{[-1^{2},0^{3};0;1^{2}]}\times 10$ &\\
\cline{2-2} &$\mathbf{[-1,0^{4};0;1^{2}]}\times 5$ &\\
\cline{2-2} &$\mathbf{[0^{7};1]} \times 1$ &\\
\cline{2-2} &$\mathbf{[0^{5};1;0;1]} \times 1$ &\\
\hline
\end{tabular}
\end{center}
\begin{itemize}
\item $\Spec(D^8_{23})=\textrm{4, 6, 7, 8, 9, 10, 12, 13, 14, 15, 16}$
\vspace{-2.5mm}\item $|Iso(D^8_{23},\cQ^8_{23})|=\textrm{960}$\vspace{-2.5mm}\item $l(D^8_{23})=\textrm{3}$
\vspace{-2.5mm}\item Antisymmetric
\vspace{-2.5mm}\item Maximally contained subpolytopes: $H(3)$, $\frac12 H(4)$, $J(7,5)$
\end{itemize}

 Perfect affine quadratic lattice $\Aff(\Z^8,\cQ^8_{
24}[x-\cc])$, where $\cQ^8_{24}$ is given by
\begin{equation*}
\setlength{\extrarowheight}{1pt} \cQ^8_{24}(\x)=\x^t\; \begin{array}{|c|c|c|c|c|c|c|c|}
\hline  24 & 18 & 18 & 18 & 18 & 18 & 27 & 13\\
\hline  18 & 24 & 18 & 18 & 18 & 18 & 27 & 13\\
\hline  18 & 18 & 24 & 18 & 18 & 18 & 27 & 15\\
\hline  18 & 18 & 18 & 24 & 18 & 18 & 27 & 19\\
\hline  18 & 18 & 18 & 18 & 24 & 18 & 27 & 15\\
\hline  18 & 18 & 18 & 18 & 18 & 24 & 27 & 15\\
\hline  27 & 27 & 27 & 27 & 27 & 27 & 39 & 21\\
\hline  13 & 13 & 15 & 19 & 15 & 15 & 21 & 19\\
\hline\end{array}
\;\x \quad \hbox{and}  \quad
\cc=\frac{1}{36}\;\begin{array}{|r|}
\hline 11\\
\hline 11\\
\hline 12\\
\hline 14\\
\hline 12\\
\hline 12\\
\hline -33\\
\hline -3\\
\hline\end{array}\;.
\end{equation*}
\begin{itemize}
\item $|O(\Z^8,\cQ^8_{24})|=\textrm{144}$; $S_{3}<O(\Z^8,\cQ^8_{24})$
\vspace{-2.5mm} \item $s(\Z^8,\cQ^8_{24})=\textrm{2}$
\vspace{-2.5mm} \item $\dim \QuadInv[\Z^8,\cQ^8_{24}]=\textrm{9}$
\end{itemize}
Inequality $\cQ_{24}^8[\x-\cc] \le \frac{97}{12}$ defines the Delaunay ellipsoid for a perfect polytope $D^8_{24}\in Del(\Z^8,\cQ^8_{24})$, whose vertex set ($|\vertex D^8_{24}|=\textrm{44}$) is given below.
\begin{center}
\setlength{\extrarowheight}{4pt}
\begin{tabular}{|c|c|c|}
\hline$x_{8}=-1$ &$x_{8}=0$ &$x_{8}=1$\\
\hline\hline$\mathbf{[0^{3};1;0^{3};-1]} \times 1$ &$\mathbf{[-1,0^{5};1;0]}\times 6$ &$\mathbf{[0^{7};1]} \times 1$\\
\cline{1-1}\cline{2-2}\cline{3-3}$\mathbf{[0^{3};1^{3};-1^{2}]} \times 1$ &$\mathbf{[0^{8}]} \times 1$ &$\mathbf{[0,1;0;-1;0^{3};1]}\times 2$\\
\cline{1-1}\cline{2-2}\cline{3-3}$\mathbf{[0^{2};1^{2};0,1;-1^{2}]}\times 2$ &$\mathbf{[0^{5},1;0^{2}]}\times 6$ &\\
\cline{1-1}\cline{2-2}$\mathbf{[0^{2};1^{4};-2;-1]} \times 1$ &$\mathbf{[0^{4},1^{2};-1;0]}\times 15$ &\\
\cline{1-1}\cline{2-2}$\mathbf{[1^{3};2;1^{2};-4;-1]} \times 1$ &$\mathbf{[0,1^{5};-3;0]}\times 6$ &\\
\cline{1-1}\cline{2-2} &$\mathbf{[1^{6};-4;0]} \times 1$ &\\
\hline
\end{tabular}
\end{center}
\begin{itemize}
\item $\Spec(D^8_{24})=\textrm{5, 8, 9, 11, 12, 13, 15, 16, 17, 19, 20, 21, 23, 24, 25, 27}$
\vspace{-2.5mm}\item $|Iso(D^8_{24},\cQ^8_{24})|=\textrm{72}$\vspace{-2.5mm}\item $l(D^8_{24})=\textrm{3}$
\vspace{-2.5mm}\item Antisymmetric
\vspace{-2.5mm}\item Maximally contained subpolytopes: $35-tope$, $H(2)$, $\frac12 H(5)$, $J(7,5)$
\end{itemize}

\include{Delaunay25_V4_GroupAttack}
 Perfect affine quadratic lattice $\Aff(\Z^8,\cQ^8_{
26}[x-\cc])$, where $\cQ^8_{26}$ is given by
\begin{equation*}
\setlength{\extrarowheight}{1pt} \cQ^8_{26}(\x)=\x^t\; \begin{array}{|c|c|c|c|c|c|c|c|}
\hline  12 & 9 & 9 & 9 & 9 & -15 & 6 & 3\\
\hline  9 & 12 & 9 & 9 & 9 & -15 & 4 & 5\\
\hline  9 & 9 & 12 & 9 & 9 & -15 & 7 & 5\\
\hline  9 & 9 & 9 & 12 & 9 & -15 & 6 & 3\\
\hline  9 & 9 & 9 & 9 & 12 & -15 & 7 & 5\\
\hline  -15 & -15 & -15 & -15 & -15 & 24 & -9 & -6\\
\hline  6 & 4 & 7 & 6 & 7 & -9 & 8 & 4\\
\hline  3 & 5 & 5 & 3 & 5 & -6 & 4 & 6\\
\hline\end{array}
\;\x \quad \hbox{and}  \quad
\cc=\frac{1}{3}\;\begin{array}{|r|}
\hline 1\\
\hline 1\\
\hline 1\\
\hline 1\\
\hline 1\\
\hline 2\\
\hline 0\\
\hline 0\\
\hline\end{array}\;.
\end{equation*}
\begin{itemize}
\item $|O(\Z^8,\cQ^8_{26})|=\textrm{2592}$; $S_{3}<O(\Z^8,\cQ^8_{26})$
\vspace{-2.5mm} \item $s(\Z^8,\cQ^8_{26})=\textrm{18}$
\vspace{-2.5mm} \item $\dim \QuadInv[\Z^8,\cQ^8_{26}]=\textrm{2}$
\end{itemize}
Inequality $\cQ_{26}^8[\x-\cc] \le 4$ defines the Delaunay ellipsoid for a perfect polytope $D^8_{26}\in Del(\Z^8,\cQ^8_{26})$, whose vertex set ($|\vertex D^8_{26}|=\textrm{45}$) is given below.
\begin{center}
\setlength{\extrarowheight}{4pt}
\begin{tabular}{|c|c|c|}
\hline$x_{8}=-1$ &$x_{8}=0$ &$x_{8}=1$\\
\hline\hline$\mathbf{[0;1;0^{2};1^{3};-1]} \times 1$ &$\mathbf{[-1,0^{4};-1;0^{2}]}\times 5$ &$\mathbf{[0;-1^{2};0;-1;-2;0;1]} \times 1$\\
\cline{1-1}\cline{2-2}\cline{3-3}$\mathbf{[0;1^{2};0^{2};1^{2};-1]} \times 1$ &$\mathbf{[0^{8}]} \times 1$ &$\mathbf{[0;-1;0^{3};-1^{2};1]} \times 1$\\
\cline{1-1}\cline{2-2}\cline{3-3}$\mathbf{[0;1^{2};0;1^{2};0;-1]} \times 1$ &$\mathbf{[0^{6};1;0]} \times 1$ &$\mathbf{[0^{7};1]} \times 1$\\
\cline{1-1}\cline{2-2}\cline{3-3}$\mathbf{[0;1^{4};2;0;-1]} \times 1$ &$\mathbf{[0^{4},1;0^{3}]}\times 5$ &$\mathbf{[1;-1;0;1;0^{2};-1;1]} \times 1$\\
\cline{1-1}\cline{2-2}\cline{3-3}$\mathbf{[0;2;0^{3};1^{2};-1]} \times 1$ &$\mathbf{[0^{3},1^{2};1;0^{2}]}\times 10$ &$\mathbf{[1;0^{2};1;0;1;0;1]} \times 1$\\
\cline{1-1}\cline{2-2}\cline{3-3}$\mathbf{[1^{3};0;1;2;0;-1]} \times 1$ &$\mathbf{[0,1^{4};2;0^{2}]}\times 5$ &$\mathbf{[1;0;1^{3};2;-1;1]} \times 1$\\
\cline{1-1}\cline{2-2}\cline{3-3} &$\mathbf{[1^{5};3;0^{2}]} \times 1$ &\\
\cline{2-2} &$\mathbf{[0^{2};1^{4};-1;0]} \times 1$ &\\
\cline{2-2} &$\mathbf{[0;1;-1;0^{3};1;0]} \times 1$ &\\
\cline{2-2} &$\mathbf{[0;1;0^{2};-1;0;1;0]} \times 1$ &\\
\cline{2-2} &$\mathbf{[1;0;1;0;1^{2};-1;0]} \times 1$ &\\
\cline{2-2} &$\mathbf{[1;0;1^{3};2;-1;0]} \times 1$ &\\
\hline
\end{tabular}
\end{center}
\begin{itemize}
\item $\Spec(D^8_{26})=\textrm{2, 3, 4, 5, 6}$
\vspace{-2.5mm}\item $|Iso(D^8_{26},\cQ^8_{26})|=\textrm{1296}$\vspace{-2.5mm}\item $l(D^8_{26})=\textrm{3}$
\vspace{-2.5mm}\item Antisymmetric
\vspace{-2.5mm}\item Maximally contained subpolytopes: $G_6$, $H(2)$, $\frac12 H(5)$, $J(7,5)$
\end{itemize}

 Perfect affine quadratic lattice $\Aff(\Z^8,\cQ^8_{
27}[x-\cc])$, where $\cQ^8_{27}$ is given by
\begin{equation*}
\setlength{\extrarowheight}{1pt} \cQ^8_{27}(\x)=\x^t\; \begin{array}{|c|c|c|c|c|c|c|c|}
\hline  12 & 9 & 9 & 9 & 9 & -15 & 7 & 3\\
\hline  9 & 12 & 9 & 9 & 9 & -15 & 6 & 2\\
\hline  9 & 9 & 12 & 9 & 9 & -15 & 6 & 3\\
\hline  9 & 9 & 9 & 12 & 9 & -15 & 6 & 3\\
\hline  9 & 9 & 9 & 9 & 12 & -15 & 8 & 5\\
\hline  -15 & -15 & -15 & -15 & -15 & 24 & -9 & -5\\
\hline  7 & 6 & 6 & 6 & 8 & -9 & 11 & 3\\
\hline  3 & 2 & 3 & 3 & 5 & -5 & 3 & 5\\
\hline\end{array}
\;\x \quad \hbox{and}  \quad
\cc=\frac{1}{12}\;\begin{array}{|r|}
\hline 1\\
\hline 3\\
\hline 2\\
\hline 2\\
\hline -2\\
\hline 1\\
\hline 3\\
\hline 3\\
\hline\end{array}\;.
\end{equation*}
\begin{itemize}
\item $|O(\Z^8,\cQ^8_{27})|=\textrm{144}$; $S_{3}<O(\Z^8,\cQ^8_{27})$
\vspace{-2.5mm} \item $s(\Z^8,\cQ^8_{27})=\textrm{24}$
\vspace{-2.5mm} \item $\dim \QuadInv[\Z^8,\cQ^8_{27}]=\textrm{8}$
\end{itemize}
Inequality $\cQ_{27}^8[\x-\cc] \le \frac{17}{4}$ defines the Delaunay ellipsoid for a perfect polytope $D^8_{27}\in Del(\Z^8,\cQ^8_{27})$, whose vertex set ($|\vertex D^8_{27}|=\textrm{44}$) is given below.
\begin{center}
\setlength{\extrarowheight}{4pt}
\begin{tabular}{|c|c|c|}
\hline$x_{8}=-1$ &$x_{8}=0$ &$x_{8}=1$\\
\hline\hline$\mathbf{[0^{4};1;0^{2};-1]} \times 1$ &$\mathbf{[-1^{5};-3;1;0]} \times 1$ &$\mathbf{[-1;0;-1^{2};-2;-3;1^{2}]} \times 1$\\
\cline{1-1}\cline{2-2}\cline{3-3} &$\mathbf{[-1,0^{4};-1;0^{2}]}\times 5$ &$\mathbf{[-1;0^{3};-1^{2};1^{2}]} \times 1$\\
\cline{2-2}\cline{3-3} &$\mathbf{[0^{8}]} \times 1$ &$\mathbf{[0^{4};-2;-1;1^{2}]} \times 1$\\
\cline{2-2}\cline{3-3} &$\mathbf{[0^{6};1;0]} \times 1$ &$\mathbf{[0^{4};-1^{2};0;1]} \times 1$\\
\cline{2-2}\cline{3-3} &$\mathbf{[0^{4},1;0^{3}]}\times 5$ &$\mathbf{[0^{7};1]} \times 1$\\
\cline{2-2}\cline{3-3} &$\mathbf{[0^{3},1^{2};1;0^{2}]}\times 10$ &$\mathbf{[0;1;0^{2};-1;0^{2};1]} \times 1$\\
\cline{2-2}\cline{3-3} &$\mathbf{[0,1^{4};2;0^{2}]}\times 5$ &$\mathbf{[0;1;0,1;0;1;0;1]}\times 2$\\
\cline{2-2}\cline{3-3} &$\mathbf{[1^{5};3;0^{2}]} \times 1$ &$\mathbf{[1^{2};0^{3};1;0;1]} \times 1$\\
\cline{2-2}\cline{3-3} &$\mathbf{[-1^{2};0^{2};-1;-2;1;0]} \times 1$ &\\
\cline{2-2} &$\mathbf{[-1;0;-1;0;-1;-2;1;0]} \times 1$ &\\
\cline{2-2} &$\mathbf{[-1;0^{2};-1^{2};-2;1;0]} \times 1$ &\\
\cline{2-2} &$\mathbf{[-1;0^{3};-1^{2};1;0]} \times 1$ &\\
\cline{2-2} &$\mathbf{[1^{4};2;3;-1;0]} \times 1$ &\\
\hline
\end{tabular}
\end{center}
\begin{itemize}
\item $\Spec(D^8_{27})=\textrm{5, 6, 7, 8, 9, 10, 11, 12, 13, 14, 15}$
\vspace{-2.5mm}\item $|Iso(D^8_{27},\cQ^8_{27})|=\textrm{72}$\vspace{-2.5mm}\item $l(D^8_{27})=\textrm{3}$
\vspace{-2.5mm}\item Antisymmetric
\vspace{-2.5mm}\item Maximally contained subpolytopes: $G_6$, $H(2)$, $\frac12 H(5)$, $J(7,5)$
\end{itemize}


\begin{thebibliography}{99}
\bibitem{Baran91}  E. P. Baranovskii (1991), Partition of Euclidean spaces into $L$-polytopes of certain perfect lattices. (Russian) Discrete geometry and
topology (Russian). \emph{Trudy Mat. Inst. Steklov. } vol. \textbf{196}, 27--46. Translated in \emph{Proc. Steklov Inst. Math. }
(1992)  vol. \textbf{196,}  no. 4, 29--51.
 
\bibitem{B} A. Barvinok (2002), \emph{A Course in Convexity,}   Graduate Studies in 
Mathematics, vol. \textbf{54}, Amer. Math. Soc., Providence, RI.

\bibitem{Barnes} E. S. Barnes (1959), The construction of  perfect and extreme forms, \emph{Acta Arithm.}, 
vol. \textbf{5}, No. 1, 57--59; No. 2, 205--222.
 
\bibitem{Atlas} J. H. Conway, R. T. Curtis, S. P. Norton, R. A. Parker, R. A. Wilson (1986),
 \emph{ATLAS of Finite Groups: Maximal Subgroups and Ordinary Characters for Simple Groups}, 
Oxford University Press. On WWW \url{http://for.mat.bham.ac.uk/atlas/} .

\bibitem{Cox1988}  H.~S.~M.~Coxeter (1988), Regular and Semiregular Polytopes III,
\emph{Math. Zeit.} vol. \textbf{200}, 3--45. Repr.\ in \emph{ Kaleidoscopes: Selected 
Writings of H.~S.~M.~Coxeter}, F.~A.~Sherk \emph{et al.}, eds., Wiley, New York, 
1995.

\bibitem{Cox1991}  H.~S.~M.~Coxeter (1991), The Evolution of Coxeter-Dynkin Diagrams,
\emph{Nieuw Archief voor Wiskunde} vol. \textbf{6}, 233--248. Repr.\ in \emph{ 
Kaleidoscopes: Selected Writings of H.~S.~M.~Coxeter}, F.~A.~Sherk \emph{et al.}, 
eds., Wiley, New York, 1995.

\bibitem{D} B.~N.~Delaunay [Delone] (1924) Sur la sph\`{e}re vide, in: \emph{Proceedings of the 
International Congress of Mathematicians,} Toronto, 1924, University of Toronto Press, Toronto, 
(1928), 695--700.

\bibitem{DGL92}  M. Deza, V. P. Grishukhin, M. Laurent (1992), Extreme hypermetrics
and $L$-polytopes. Sets, graphs and numbers (Budapest, 1991),
157--209, \emph{Colloq. Math. Soc. Janos Bolyai}, vol. \textbf{60},
North-Holland, Amsterdam.

\bibitem{DLbook}  M. Deza and M. Laurent (1997), \textit{Geometry of cuts and metrics},
Algorithms and Combinatorics \textbf{15}, Springer-Verlag, Berlin.

\bibitem{D05}  M. Dutour (2005), Infinite serie of extreme Delaunay polytopes, \emph{European Journal of Combinatorics }
vol. \textbf{26,} no. 1, pp. 129--132.

\bibitem{D04}  M. Dutour (2004), The six-dimensional Delaunay polytopes,
\emph{Europ. J. Comb.} vol. \textbf{25 }, 535--548.

\bibitem{Dprogram}  M. Dutour (2002), \emph{EXT-HYP7}, Program proving the uniqueness of 
perfect Delaunay polytope in $\R^6$, \url{http://www.liga.ens.fr/~dutour/HYP7/index.html}

\item M. Dutour and F. Vallentin (2005), Some six-dimensional rigid lattices, \emph{Voronoi's 
Impact on Modern Science, Book 3. Proc. Inst. Math. Nat. Acad. Sci. Ukraine} (2005), vol.
\textbf{55}, 102--108.

\bibitem{DSV}  M. Dutour, A. Schuermann, F. Vallentin, (2006)  \emph{A generalization of Voronoi's reduction theory and its application}, 
Preprint \url{http://arXiv.org/math.MG/0601084}.

\bibitem{E75}  R. Erdahl (1975), \emph{A convex set of second-order inhomogeneous
polynomials with applications to quantum mechanical many body theory, } Mathematical 
Preprint \#1975-40, Queen's University, Kingston, Ontario.

\bibitem{E92}  R. Erdahl (1992), A cone of inhomogeneous second-order polynomials,
\textit{Discrete Comput. Geom}. vol. \textbf{8}, no. 4, 387--416.


\bibitem{EOR_05} R.M. Erdahl, A. Ordine, K. Rybnikov (2005), Perfect Delaunay Polytopes, 
\emph{Voronoi's Impact on Modern Science, Book 3. Proc. Inst. Math. Nat. Acad. Sci. Ukraine} 
(2005), vol. \textbf{55}, 126--136.

\bibitem{EOR_04} R.M. Erdahl, A. Ordine, K. Rybnikov (2004), \emph{Constructions for Perfect 
Quadratic Functions and Delaunay Polytopes,} submitted, \url{http://arxiv.org/math.NT/0408122} 
.

\bibitem{ER01-02}  R. M. Erdahl and K. Rybnikov (2002), \emph{Supertopes,}  
\url{http://arXiv.org/math.NT/0501245} .

\bibitem{ER02}  R. M. Erdahl and K. Rybnikov (2002), Voronoi-Dickson Hypothesis on Perfect 
Forms and L-types, Peter Gruber Festshrift: \emph{Rendiconti del Circolo Matematiko di 
Palermo,} Serie II, Tomo \textbf{LII,} part I, 279--296. 

\bibitem{Gosset} T.~Gosset (1900), On the regular and semi-regular figures in space of $n$ dimensions,
\emph{Messenger of Math.} vol. \textbf{29}, 43--48.

\bibitem{} P. Gruber and C. Lekkerkerker  (1987), \textit{Geometry of Numbers,} 2-nd edition, Elsevier Science
Publishers.

\bibitem{Humphreys}  J. E. Humphreys (1990), \emph{Reflection Groups and Coxeter Groups,} 
Cambridge Studies in Advanced Mathematics, vol. \textbf{29}, Cambridge University Press, 
Cambridge.

\bibitem{KZ}  A. Korkine and G. Zolotareff (1873), Sur les formes quadratiques,
\textit{Math. Ann}. no. 6, 366--389.

\bibitem{M03}  J. Martinet (2003), \textit{Perfect lattices in Euclidean spaces},
Fundamental Principles of Mathematical Sciences \textbf{327}, Springer-Verlag, 
Berlin.

\bibitem{ER01}  K. Rybnikov (2001), \emph{REU 2001 Report: Geometry of Numbers}, Research Experience for Undergraduates
 Project, Final Report, Department of Mathematics,  Cornell University
\url{ http://www.mathlab.cornell.edu/~upsilon/REU2001.pdf} .

\bibitem{ER01}  S. S. Ryshkov and S. Sh. Shushbaev (1981), The structure of the L-partition for the second perfect lattice (Russian), 
\emph{Matematicheskii Sbornik} (Russian), vol. \textbf{116,} No. 2, 218--231.
\end{thebibliography}
\end{document}